\documentclass[onefignum,onetabnum]{siamart171218}

\usepackage{lipsum}
\usepackage{amsfonts}
\usepackage{amssymb}
\usepackage{graphicx}
\usepackage{epstopdf}
\usepackage{algorithmic}
\ifpdf
  \DeclareGraphicsExtensions{.eps,.pdf,.png,.jpg}
\else
  \DeclareGraphicsExtensions{.eps}
\fi

\usepackage{subcaption}

\newsiamremark{remark}{Remark}
\newsiamremark{hypothesis}{Hypothesis}
\crefname{hypothesis}{Hypothesis}{Hypotheses}
\newsiamthm{claim}{Claim}

\headers{Optimal control of COVID-19 infection rate with social costs}{Aaron Z.\ Palmer,  Zelda B.\  Zabinsky, and Shan Liu}

\title{Optimal control of COVID-19 infection rate considering social costs \thanks{This work has been funded in part by the National Science Foundation grant CMMI-1935403.}}

\author{Aaron Z.\ Palmer\thanks{Mathematics Department, University of British Columbia, Vancouver BC, Canada
  (\email{azp@math.ubc.ca}).}
\and Zelda B. Zabinsky\thanks{Department of Industrial \& Systems Engineering, University of Washington, Seattle WA
  (\email{zelda@uw.edu}, \email{liushan@uw.edu}).}
\and Shan Liu\footnotemark[3]}

\usepackage{amsopn}

\ifpdf
\hypersetup{
   pdftitle={Optimal control of COVID-19 infection rate considering social costs},
  pdfauthor={Aaron Z.\ Palmer,  Zelda B.\ Zabinsky, and Shan Liu}
}
\fi

\newcommand{\R}{\mathbb{R}}

\newcommand{\ba}{\begin{aligned}}
\newcommand{\ea}{\end{aligned}}
\newcommand{\bt}{\begin{thm}}
\newcommand{\et}{\end{thm}}
\newcommand{\bc}{\begin{corollary}}
\newcommand{\ec}{\end{corollary}}
\newcommand{\bl}{\begin{lemma}}
\newcommand{\el}{\end{lemma}}
\newcommand{\bpf}{\begin{proof}}
\newcommand{\epf}{\end{proof}}
\newcommand{\bpb}{\begin{problem}}
\newcommand{\epb}{\end{problem}}
\newcommand{\bd}{\begin{definition}}
\newcommand{\ed}{\end{definition}}
\newcommand{\bn}{\begin{note}}
\newcommand{\en}{\end{note}}
\newcommand{\bq}{\begin{question}}
\newcommand{\eq}{\end{question}}
\newcommand{\bp}{\begin{proposition}}
\newcommand{\ep}{\end{proposition}}
\newcommand{\be}{\begin{example}}
\newcommand{\ee}{\end{example}}
\newcommand{\bex}{\begin{exercise}}
\newcommand{\eex}{\end{exercise}}
\newcommand{\ben}{\begin{enumerate}}
\newcommand{\een}{\end{enumerate}}
\newcommand{\nn}{\nonumber}

\begin{document}

\maketitle

% REQUIRED
\begin{abstract}
    The COVID-19 pandemic has posed a policy making crisis where efforts to slow down or end the pandemic conflict with economic priorities.  This paper provides mathematical analysis of optimal disease control policies with idealized compartmental models for disease propagation and simplistic models of social and economic costs.  Two locally optimal control strategies are found and categorized as `suppression' and `mitigation' strategies. We analyze how these strategies change when we incorporate vaccination into the model and find a new optimal `delay-mitigation' strategy.
\end{abstract}

% REQUIRED

% REQUIRED

% {\red 
% To DO:
% \begin{enumerate}

%   \item Revise Introduction/Conclusion/Abstract
%   \item Replace Simplified SIR figure with Vaccination figure.
%   \item Revise New Vaccinations Section
%   \item Consider the parameter $\delta_0$ and justify (or change and justify it)
%   \item Possibly revise the appendix sections

% \end{enumerate}
% }

\section{Introduction}

The spread of infectious diseases in a population can be modeled in a variety of ways.  Here we focus on the optimization and control of stochastic and deterministic `compartmental' models, where the population is grouped into a handful of states that represent the progression of the disease.  These models have the benefit of having relatively few variables and parameters, and can be comprehended intuitively.  

Our main contribution is to study the optimal control of models pertaining to COVID-19 with control implemented through influencing the infection rate, e.g., by social distance policies. For instance, stay-at-home orders, wearing masks, and other non-pharmaceutical interventions (NPIs) can serve to reduce the infection rate, while relaxing NPIs and opening up society increases infection rate. Our analysis focuses on qualitative aspects of the optimal control policies with an aim to better understand the effective strategies for combating COVID-19.  {In the optimal control model, we optimize over a combination of: (i) the end time, which represents the time the disease is below a given threshold and considered under control; (ii) a running cost due to controlling infection rate, that may be interpreted as the societal cost of implementing NPIs; (iii) a running hospitalization cost; and (iv) a death cost at end time.

We use realistic parameters for COVID-19, and}  solve for optimal control policies. The optimal polices fall into one of two categories that are widely used to describe real world responses: `suppression' strategies that aim to eliminate the disease in a population, and `mitigation' strategies that aim to reduce the negative impact of the disease. {We  find that the `suppression' strategy performs significantly better than the `mitigation' or any other strategy.  

We then incorporate vaccinations and find that while the `suppression' strategy remains a viable solution under some conditions, the `mitigation' strategy transforms into a globally optimal `delay-mitigation' strategy, where the spread of the disease is delayed until vaccination takes effect.}

\subsection{Literature Review}
There is a long history of `compartmental' models including previous research into optimization and control and applications to the COVID-19 pandemic.  We provide a non-exhaustive summary of relevant works.  {In these works it can be seen how different specifications of the problem yield distinct results.  Despite the differences, these results share important qualitative features that are also reflected in the results of our approach.} 

{A classic optimal control formulation for} deterministic and stochastic Susceptible-Infectious-Recovered (SIR) compartmental models is analyzed in \cite{wickwire1975optimal}.  The control variable adjusts the rate at which to remove individuals from the infected population by isolation.  The solutions that are found either expend all or none of the resources. 
Deterministic SIR and Susceptible-Exposed-Infectious-Recovered (SEIR) models are considered in \cite{behncke2000optimal} with control by vaccination, quarantine, screening, or health campaigns.  It is found that in all cases the optimal policy to minimize number of infected is to apply maximum effort on an initial time interval.  {A similar finding was made for an SEIR type model applied to COVID-19 in \cite{perkins2020optimal} with control on the {infection}  rate}. An SEIR model with logistic population growth is considered in \cite{thater2018optimal}, where the control found by numerical optimization appears to approach a state where the disease is endemic in the population.
An SEIR type model has been applied to COVID-19 in \cite{djidjou2020optimal} with the infection rate controlled up until a vaccine is developed. The optimal control found suppresses the disease until near the development of the vaccine when the control is relaxed. A variation is considered \cite{ketcheson2020optimal} where instead of vaccination it is assumed that control reverts back to the norm, and the result is a mitigation strategy that brings the susceptible population to herd immunity. 
A deterministic SIR model with control by vaccination and isolation is considered in \cite{hansen2011optimal}.  They minimize total outbreak size, and show that it is optimal to allocate the available resources to until the resources run out. 

A discounted deterministic Susceptible-Infectious-Susceptible (SIS) model with two {geographic regions appears in \cite{rowthorn2009optimal}, where the allocation of medical resources between the two regions is optimized.}  The main finding is that {it is better to treat the region with less {incidence of} disease first.} 
A review of epidemic models and control was undertaken in \cite{nowzari2016analysis} with an emphasis of effects from the network of individuals (not seen in compartmental models). 
A spatial SIR model is considered in \cite{lee2020controlling} with control of spatial dynamics. 

The distinction between open loop and feedback control for a stochastic epidemic simulation model is studied in \cite{bussell2019applying}, where the {feedback control} is found to perform better because the approximate model loses accuracy across the entire time domain. {In this paper, we focus on open loop control and consider some aspects of feedback control in the appendix.}

The report \cite{ferguson2020report} provides an analysis of `suppression' and `mitigation' strategies for the COVID-19 pandemic, although optimality is not addressed.  Numerous papers have analyzed the efficacy of different non-pharmaceutical interventions for COVID-19, including  \cite{acuna2020modeling}, \cite{li2020modeling}, \cite{ngonghala2020mathematical}, \cite{tuite2020mathematical}.   We do not address specific non-pharmaceutical interventions, but recognize that various interventions (e.g., social distancing, mask wearing, lock-down orders) impact the infection rate, which is our control variable.
%For  example the finding of \cite{ngonghala2020mathematical} that population-wide use of low efficacy  masks would slow but not stop the spread of COVID-19, whereas more costly high efficacy masks may lead to elimination of COVID-19 if used widely.  
The recent work \cite{naraigh2020piecewise} considers optimal control using piecewise constant strategies, and compares the optimal `suppression' strategy with a non-optimal `mitigation' strategy.  In our results, both the `suppression' strategy and the `mitigation' strategy are found to be locally optimal, however, the end time for the `mitigation' strategy is much longer than that of the `suppression' strategy. 

We summarize the elements of our model and results that have not appeared in prior works:
\begin{itemize}
  \item The form of the cost function includes a social and economic control cost that represents the impact of NPI's, as well as hospitalization and death costs. The optimal control policy maintains positivity of the controlled infection rate while allowing for a continuum of values.
    \item The end time in our model is determined by elimination of the disease (below a threshold) in the population.
  \item We find two locally optimal strategies, both of which are the global optimum for different choices of parameters.
  \item We combine the control of infection rate with a vaccination roll-out to find a qualitatively novel strategy.
\end{itemize}

\subsection{Model Assumptions}
We consider a model that has six states: susceptible, exposed, infectious, hospitalized, recovered, and dead.  We optimize a cost comprised of the cost of hospitalizations, the cost of deaths, and the social/economic costs of reducing the infection rate by social distancing policies and other NPIs.  The model terminates when the disease has {been eliminated from} the population (below a small threshold).     We also incorporate a vaccine, that removes individuals from the susceptible population at a fixed rate.   We assume that the entire susceptible population is willing to get vaccinated, and that, once vaccinated, achieve immunity.

Realistic parameters for transition rates and costs are chosen from related literature and data on COVID-19.  We also test the model with ranges of values for some parameters, including population size, mean infectious period, cost coefficient, vaccination rate, and initial conditions.  
%For the results without vaccination, we consider the population size of Washington State and use initial conditions corresponding to COVID-19 levels in Washington State in mid 2020.  For the results with vaccination, we scale up the population size to California State, and to the population size of the USA, with initial conditions corresponding to COVID-19 levels in the USA in early 2021.

Our model assumes a negligible portion of the population can be re-infected with the disease after recovering, and omits births and immigration.  By assuming that the parameters of the disease remain constant, we ignore the possibility that mutations of the virus may change its parameters or allow the mutated virus to infect previously recovered or vaccinated individuals.  In the appendix, we propose a stochastic optimal control model, which could incorporate uncertainty due to mutations as well as uncertainty in the development of future vaccines.

\subsection{Results}
Our first result is to find two locally optimal strategies in the deterministic model that we characterize as `suppression' and `mitigation' strategies. The locally optimal strategies balance the cost of the control (e.g., social and economic cost) with the cost of hospitalization and death. The  `suppression' strategy reduces the reproduction number well below 1 and it remains below 1 until the disease has been eliminated.   This strategy is qualitatively similar to the strategies in the literature that apply maximum control at the beginning of an outbreak, and also resembles the strategy up to vaccination time found in \cite{djidjou2020optimal}.    The `mitigation' strategy instead applies control around the peak of the epidemic rather than at the initial time. The reproduction number is then brought below 1 by a combination of control and the development of herd immunity. By the time the disease is eliminated from the population, the majority of the population will have gotten infected when following the `mitigation' strategy. This strategy is similar to what is found in \cite{ketcheson2020optimal}.  We find both of these solutions exist as local optima with reasonable parameters for COVID-19, while the suppression strategy is the global minimum.

When accounting for the roll-out of vaccinations, we find that the `suppression' strategy remains a local optimum, as long as the time to elimination is shorter than the roll-out time of the vaccine. However, a new strategy that we call `delay-mitigation' is the global optimum.  The `delay-mitigation' is a more complex form of the `mitigation' strategy that combines early suppression with later mitigation as the vaccine rolls out.  In this new strategy, the reproduction number is maintained near 1 until vaccination of the population drives it down. The delay-mitigation strategy is consistent with U.S. public health messaging that the control (e.g., NPIs) can be relaxed when a large portion of the population is vaccinated.

We also consider a stochastic model in the appendix that naturally provides a globally optimal feedback control and yields threshold values of infected populations for which to apply different policies.  For example, when the infectious population reaches a threshold,  more control is exercised (in the form of NPIs) at a social and economic cost. This model captures the fluctuations of the end time, which leads to high variability of the cost for the `suppression' strategy. 
The development of a vaccine allows the optimal feedback control to lessen, however, when the number of cases is still small, the optimal feedback control strengthens to reduce the infection rate upon development of the vaccine.

The code (written in Python with NumPy and Matplotlib for figures) and numerically generated data (JSON format) that we used is available at \url{https://github.com/AaronZPalmer/SEIHRD.git}.

\section{Epidemic Model}

\subsection{Susceptible-Exposed-Infectious-Hospitalized-Recovered-Dead \\ (SEIHRD) Model}

We assume that individuals are indistinguishable and reside in one of six states: susceptible, exposed, infectious, hospitalized, recovered, and dead.

The state variables, which represent number of individuals in each state at a given time, are: 
\begin{itemize}
    \item $S_t$, the number of `susceptible' individuals at time $t$;
    \item $E_t$, the number of `exposed' individuals at time $t$;
    \item $I_t$, the number of `infectious' individuals at time $t$;
    \item $H_t$, the number of `hospitalized' individuals at time $t$;
    \item $R_t$, the number of `recovered' individuals at time $t$;
    \item $D_t$, the number of `dead' individuals at time $t$.
\end{itemize}
The total population size $N$ accounts for all individuals $N = S_t + E_t + I_t + H_t + R_t + D_t$.

We assume that the following individual transitions occur at exponentially distributed times:
\begin{itemize}
    \item Susceptible to exposed at rate $\beta\, \frac{I_t}{N}$;
    \item Exposed to infectious at rate $\alpha$;
    \item Infectious to hospitalized at rate $\lambda_0$;
    \item Infectious to recovered at rate $\gamma_0$
    \item Infectious to dead at rate $\delta_0$;
    \item Hospitalized to recovered at rate $\gamma_1$;
    \item Hospitalized to dead at rate $\delta_1$.
\end{itemize}

\begin{figure}
\centering
\includegraphics[width=0.4\textwidth]{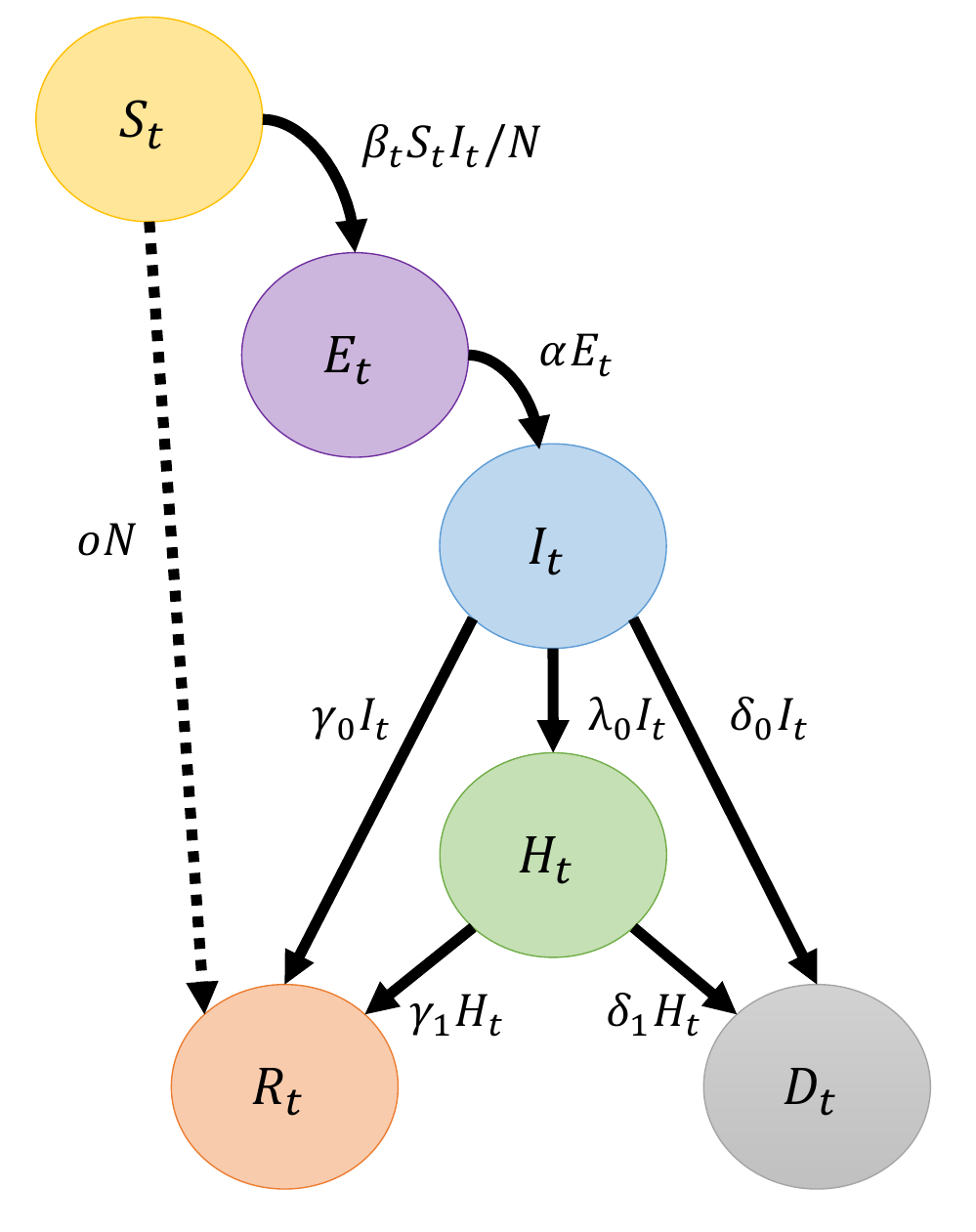}
\caption{Rate transition diagram for SEIHRD model. The dashed line represents the rate transition as vaccines roll-out.}
\label{fig:ratediagram}
\end{figure}

The transition of any susceptible person to become infected, in other words $S_t$ becomes $S_t - 1$ occurs as the first of $S_t$ independent random times with exponential distributions of rate 
$\beta\, \frac{I_t}{N}$.  
The first of these times that marks the transition from $S_t$ to $S_t-1$ is exponentially distributed with rate 
$\beta\, S_t\, \frac{I_t}{N}$.  
The same holds for the other transitions.  An illustration of the SEIHRD  rate transition diagram is in Figure~\ref{fig:ratediagram}. The dashed line in Figure~\ref{fig:ratediagram} represents the transition due to vaccination, discussed in Section~\ref{Sec:vaccinations}.

When $N$ is large and the macroscopic state variables are of order $N$, we can approximate the epidemic dynamics by the system of differential equations
\begin{align}\label{eqn:mf_de}
    \frac{dS_t}{dt} =& - \beta_t\, S_t\, \frac{I_t}{N}\\
    \frac{dE_t}{dt} =&\ \beta_t\, S_t\, \frac{I_t}{N} - \alpha\, E_t\nn\\
    \frac{dI_t}{dt} =&\ \alpha\, E_t - \gamma_0\, I_t-\lambda_0\, I_t-\, \delta_0\, I_t\nn\\
    \frac{dH_t}{dt} =&\ \lambda_0\, I_t - \gamma_1\, H_t- \delta_1\, H_t\nn\\
    \frac{dR_t}{dt} =&\ \gamma_0\, I_t + \gamma_1\, H_t\nn\\
    \frac{dD_t}{dt} =&\ \delta_0\, I_t + \delta_1\, H_t.\nn
\end{align}

A common theme of epidemic models is that the dynamic behavior is characterized by a single `reproduction number,' which when greater than 1 corresponds to exponential growth and when less than 1 corresponds to exponential decline of the disease in the population. 
The reproduction number ${\bf R_0}$ is the average number of infections caused by a single infected individual in a susceptible population.  This is calculated, in our notation, as the infection rate $\beta$ times the expected duration of the disease, 
\begin{align}\label{eqn:R_0}
    {\bf R_0} = \frac{\beta}{\lambda_0+\delta_0+\gamma_0}.
\end{align} 
The effective reproduction number, in our notation, is
\begin{align}\label{eqn:R_e}
    {\bf R_e} = {\bf R_0}\, \frac{S}{N}= \frac{\beta\, \frac{S}{N}}{\lambda_0+\delta_0+\gamma_0},
\end{align}  
{that is, the number of infections generated in the current state of a population (rate of infection per infected individual, $\beta\, S / N$, times average infection duration $1/(\lambda_0+\delta_0+\gamma_0)$),   under the simplifying assumption that $S$ is constant in the relevant time interval.} The effective reproduction number governs whether the disease level in the population will increase, ${\bf R_e}>1$, or decrease, ${\bf R_e}<1$.

\subsection{Equilibria}

After a sufficiently long time the solutions of (\ref{eqn:mf_de}) will approach equilibrium points where the right hand sides of the differential equations are zero.  An equilibrium requires that $E = I = H = 0$, since $\frac{dH}{dt}=0$ and $\frac{dR}{dt}=0$ imply $I=H=0$, and then $\frac{dE}{dt}=0$ implies $E=0$. It is then clear that $S$, $R$, and $D$ can take any values in equilibrium and this allows us to determine all of the equilibrium points.  In Appendix \ref{apx:equilibria} we further investigate the stability of the equilibria following standard techniques and provide illustrations of how the SEIHRD model reaches an equilibrium.

\subsection{Parameter Settings for the SEIHRD Model of COVID-19}

We use published parameters from scientific papers and in some cases, parameters based on COVID-19 data for the United States \cite{CDC-burden}.  
%(see \url{https://www.cdc.gov/coronavirus/2019-ncov/cases-updates/burden.html}. ). 
Each parameter in Table~\ref{table:parameters} is a rate per day per person except for $N$. 
Some of the parameters can be obtained directly from the literature and other parameters are calculated. We also perform a sensitivity analysis on several parameters.

To determine $\delta_1$, we use  $\gamma_1=0.1$ and the value 11.8\% for the percentage of hospitalizations that lead to death (from \cite{moghadas2020projecting})  and set the ratio $\frac{\delta_1}{\gamma_1+\delta_1}=0.118$. 
Similarly, we note that \cite{moghadas2020projecting} provides $\lambda_0+\gamma_0+\delta_0=0.217$ and  \cite{CDC-burden} provides the percentage of infections that lead to hospitalization as 3.69\%, then we use the   ratio $\frac{\lambda_0}{\lambda_0+\gamma_0+\delta_0}=0.0369$ to calculate $\lambda_0$.

 To determine $\delta_0$ we calculate the percentage of infections that lead to death by
$$	\frac{\delta_0 + \lambda_0 \left(\frac{\delta_1}{\gamma_1+\delta_1}\right)}{\lambda_0 + \gamma_0 + \delta_0} = \frac{\rm Total\ Deaths}{\rm Total\ Cases} = \frac{1.4\, *\, 346,927}{4.6\, * \, 20,112,544 }
$$
where total deaths and total cases are reported numbers from the United States as of January 1, 2021 \cite{CDC-tracker}, and 40\% additional deaths is a  reported number in \cite{cutler2020covid}, and 4.6 times as many cases than confirmed cases comes from  \cite{CDC-burden}.

We treat the infection rate $\beta$ as our control variable, and note that
a range of infection rates between $0.11$ and $0.87$ corresponds to 
values of   ${\bf R_0}$ from $0.5$ to~$4$.

\begin{table}[ht]
\caption{Parameter settings}
\label{table:parameters}
\centering
\begin{tabular}{| l | l | l |}
\hline
$\alpha$ & $0.192$ & The exposed to infectious rate is the reciprocal of \\
& & the mean incubation period of $5.2$ days. \cite{moghadas2020projecting} \\
\hline
%$\beta$ & $0.11 \sim 0.87$ & A range of infection rates, corresponding to \\
%& &  values of   ${\bf R_0}$ in the range $0.5 \sim 4$. \cite{CDC-R0, sanche2020-R0} \\
%\hline
$\lambda_0+\gamma_0+\delta_0$ & 0.217 & The total rate of leaving state $I$  is the reciprocal   \\
& & of the mean infectious period, 4.6 days. \cite{moghadas2020projecting} \\
& & We also consider a range of $4$ to $10$ days.\\
\hline
$\lambda_0$ & 0.008 & The hospitalization rate is calculated from \\
& &  3.69\% hospitalized.  \cite{CDC-burden} \\
\hline
$\delta_0$ & 0.000195& The death rate is calculated from US data, \cite{CDC-tracker} \\
& &  and assumes 40\% excess deaths. \cite{cutler2020covid} \\
\hline
$\gamma_0$ & 0.209 & The recovery rate is the remainder, such that\\
& &   $\lambda_0+\gamma_0+\delta_0$ equals total rate of leaving state $I$.\\
\hline
$\gamma_1$ & 0.1 & The recovery rate from the hospital is the reciprocal \\
& & of the mean length of hospital stay, 10 days. \cite{moghadas2020projecting} \\
\hline
$\delta_1$ & 0.013 &  The hospital death rate is calculated from  11.8\%\\
& & of those hospitalized die, so that $0.118 = \frac{\delta_1}{\gamma_1+\delta_1}$. \cite{moghadas2020projecting}\\
\hline
$N$ & 7,600,000& Population of Washington State. \\
&& We also consider a range from $1$ million to \\  
& & the U.S. population  ($328.2$ million), \\
& & and the world population ($7.8$ billion).\\
\hline
\end{tabular}
\end{table}

We choose the initial values consistent with COVID-19 infection levels in Washington State, as of June 1, 2020. 
%consistent with the parameters above, and 
% with 0.02\% infected corresponding to COVID-19 levels in Washington State in mid 2020.  
We set  $D_0=0$ in the initial conditions, so  $D_t$ represents cumulative  deaths since June 1, 2020 in our criterion.
The initial conditions for Washington State, based on June 1, 2020 are 
\begin{align}\label{eqn:initial_values}
	\left(\begin{array}{c}
		S_0\\
		E_0\\
		I_0\\
		H_0\\
		R_0\\
		D_0
	\end{array}\right) = 
	\left(\begin{array}{c}
		7,600,000 - E_0 - I_0 - H_0 - R_0\\
		294 * 4.6 / \alpha\\
		294 * 4.6 / (\lambda_0+\gamma_0+\delta_0)\\
		178 * 1.9\\
		22,238 * 4.6 - E_0 - I_0 - H_0\\
		0
	\end{array}\right)
	=\left(\begin{array}{c}
		7,497,705\\
		7,044\\
		6,221\\
		338\\
		88,692\\
		0
	\end{array}\right).
\end{align}
We  used 294 for daily confirmed cases from the 7-day rolling average of daily confirmed cases on June 1, 2020, 178 current hospitalizations as of May 31, 2020, and 22,238 total cases as of June 1, 2020 in Washington State \cite{WA-health}. We  used the estimate of 4.6 times more cases than confirmed cases and 1.9 times more hospitalized cases than reported, as in \cite{CDC-burden}.
% {\red When the population size is scaled up, the initial conditions are scaled up proportionally.}
%\begin{align}
%%\label{eqn:initial_values}
%\left(\begin{array}{c}S_0 \\ E_0 \\I_0\\ H_0\\ R_0\\ D_0\end{array}\right)=\left(\begin{array}{c} 7,596,654 \\ 1,718 \\1,520\\ 108\\ 0\\ 0\end{array}\right).
%\end{align}
These initial conditions are used in the examples unless otherwise stated.

\section{Optimal Control} \label{sec:det}

We now pose an optimal control problem supposing that the infection rate $\beta_t$ can be controlled  over time through social distancing and other non-pharmaceutical interventions.

We assume the dynamics of (\ref{eqn:mf_de}), and we optimize over 
the end time $T$ and control policy $(\beta_t)_{t\in [0,T]}$. The cost consists of a running control cost $L$, a running hospitalization cost $F$, and a terminal death cost $G$, which has the form 
\begin{align}\label{eqn:cost}
    J\big[(\beta_t)_{t\in [0,T]}, T\big] = \int_0^T\Big[L(\beta_t)+F(H_t)\Big]dt + G(D_T).
\end{align}
We constrain the final state by 
\begin{align}\label{eqn:extinction}
    E_T + I_T + H_T \leq e^{-1}
\end{align}
 to represent that the disease has {been eliminated from the population}. The end time $T$ is the first time for which this constraint is satisfied.  {The threshold value of $e^{-1}$ is used in (\ref{eqn:extinction}) so that the end time $T$ coincides with the expected time that $E_t=I_t=H_t=0$ in a probabilistic model where the end time is exponentially distributed.  For example if $I_0=1$, then the expected time to transition to $I_t=0$ coincides with the time $t$ such that  $\mathbb{E}[I_t]=e^{-1}$.}  
% {We note that the term extinction time in  public health terminology is only  used when dealing with the entire global population}. 

  Thus the optimal control problem has the following elements:
\begin{itemize}
    \item Decision variables $T$ and $(\beta_t)_{t\in [0,T]}$;
    \item Dynamics determined by (\ref{eqn:mf_de});
    \item Objective function (\ref{eqn:cost});
    \item Terminal constraint (\ref{eqn:extinction}).
\end{itemize}

Note that while the solution to the SEIHRD model scales linearly with the population size $N$, the {end} time is nonlinear in $N$. This is important in interpreting the optimal control policy, and its dependence on population size and end time.  {For this reason, and to better illustrate the relative scale of the solution, we have left the state variables as total number of individuals rather than rescaling to dimensionless ratios like $S_t/N$, etc.}

We select cost functions for the optimal control criterion in (\ref{eqn:cost}) by a phenomenological approach:
\begin{align}
L(\beta)=&\ N\, k\, \Big(-\log\big(\frac{\beta}{b}\big) + \frac{\beta}{b} - 1\Big), \label{eqn:costL}\\
F(H)=&\ c_0\, H+\frac{c_1}{N}\, {H}^2, \label{eqn:costF}\\
G(D)=&\ d\, D. \label{eqn:costG}
\end{align}
%{ We choose a convex form for the control cost $L(\beta)\geq 0$ when $\beta \geq 0$, with $L(b)=0$, and  
% assume that $L(\beta)=+\infty$ if $\beta\leq 0$.}

 {The control cost is a function of the control $\beta$,  the baseline infection rate $b$, and a control cost coefficient $k$ that reflects the cost of a proportional reduction in ${\bf R_0}$.
 We choose a convex form for the control cost $L(\beta)\geq 0$ when $\beta \geq 0$ and  
 assume that $L(\beta)=+\infty$ if $\beta\leq 0$. We set $L(b)=0$, where the parameter $b$ is the baseline infection rate, so there is no additional control cost.}
 {The logarithmic term in $L$ reflects that the cost to reduce the infection rate (and hence the reproduction factor) by 1\% is a constant $ N\, k /100$.  The linear and constant terms of $L$ normalize so that the minimum is at $\beta=b$ and has zero cost; values of $L$ for $\beta\geq b$ are not important as it is never optimal to increase the infection rate beyond the natural baseline infection rate.  This choice of a logarithmic term in $L$ is consistent with the phenomenon that independent methods of intervention have additive costs and multiplicative reduction of $\beta$. A more common quadratic optimal control cost, e.g., $(\beta-b)^2$, should not be used here as it fails to enforce  the infection rate to remain positive, unrealistically allowing the infection rate to be zero with a finite cost.} 
 
 A quadratic term is included in the hospitalization cost $F$ to reflect the cost of passing a hospital occupancy threshold.  The linear death cost reflects each human life being of equal value.

\subsection{Control Parameter Settings}

The cost functions of (\ref{eqn:costL}), (\ref{eqn:costF}) and (\ref{eqn:costG}) have {five} parameters: the control cost coefficient $k$, the baseline infection rate $b$, the hospitalization cost rates {$c_0,c_1$} and the death cost coefficient $d$.  We attempt to choose reasonable values in USD, shown in Table~\ref{table:control_parameters}, though further specification is ultimately subjective. {The control cost coefficient  $k$ reflects the cost of a proportional reduction in ${\bf R_0}$. We tested a range of values for $k$,  from $50$ to $250$.}  It is hopeful that an efficient handling of the epidemic would significantly lower the control cost coefficient, $k$. We also vary the population size $N$ and note the resulting end time.

\begin{table}[ht]
\caption{Control parameter settings}
\label{table:control_parameters}
\centering
\begin{tabular}{| l | l | l |}
%\hline
%$k$ & 100 & An approximate value is chosen such that, at\\
%&We also consider &${\bf R_0} =1/2$, the cost is \$120 a day per person \\
%&50 $\sim$ 250 & ${\bf R_0}=1$, the cost is \$64 a day per person,\\
%&  & consistent with analysis in \cite{cutler2020covid}. \\
%&&An approximate value is chosen such that, at\\
%$k$ &50 \sim 250, where \$100 per day per person&  A subjective value chosen to be \\ 
%& (roughly \$1 per day per person & a round ballpark number; \\
%& for a 1\% reduction of $R_0$) & at $R_0=1/2$, the cost is \$120 a day per person;\\
%& & at $R_0=1$, the cost is \$64 a day per person.\\
\hline
$b$ & $0.87$& Choose $b=\beta$ such that ${\bf R_0}= 4$.\\
\hline
$k$ & 100 & An approximate value is chosen such that, when\\
& &${\bf R_0} =1/2$, the cost $L(\beta)= \$120 N$  \\
&  & and when ${\bf R_0}=1$, the cost is $L(\beta)= \$64 N$,\\
&  & consistent with analysis in \cite{cutler2020covid}. \\
\hline
$c_0$ & 3,500 & Calculated for a total cost per patient of \\
& & \$35,000 divided by an average stay \\
& & of 10 days.  \cite{Rae-hospcosts}  \\
\hline
$c_1$ & \ $c_0/2$  & Chosen to be the same order of magnitude \\
& &  as $c_0$. \\
\hline
$d$ & 7,000,000 & The cost of a single death used in \cite{cutler2020covid}.  \\
\hline
\end{tabular}
\end{table}

%\begin{table}[ht]
%\caption{Control parameter settings}
%\label{table:control_parameters}
%\centering
%\begin{tabular}{| l | l | l |}
%\hline
%$d$ & \$1,000,000 per person& US department of transportation values   \\
%& &  life at \$9,600,600; this accounts for  \\
%& &approximately 10\% shortened lifespan.\\
%\hline
%$b$ & $0.87$ per day per person& For $R_0= 4$.\\
%\hline
%$c_0$ & \$3,500 per day per person & Calculated for a total cost per patient of \\
%& & \$35,000 divided by an average stay \\
%& & of 10 days.  \cite{Rae-hospcosts}  \\
%\hline
%$c_1$ & \ $c_0/2$  & Chosen to be the same order of magnitude \\
%& &  as $c_0$. \\
%\hline
%$k$ & 50 $\sim$ 250, where  & A subjective value chosen to be \\
%&\$100 per day per person is &a round ballpark number; at\\
%& roughly \$1 per day per person  & $R_0=1/2$, the cost is \$120 a day per person\\
%&for a 1\% reduction of $R_0$ &$R_0=1$, the cost is \$64 a day per person.\\
%%$k$ &50 \sim 250, where \$100 per day per person&  A subjective value chosen to be \\ 
%%& (roughly \$1 per day per person & a round ballpark number; \\
%%& for a 1\% reduction of $R_0$) & at $R_0=1/2$, the cost is \$120 a day per person;\\
%%& & at $R_0=1$, the cost is \$64 a day per person.\\
%\hline
%\end{tabular}
%\end{table}

\subsection{Deterministic Dynamics Results}\label{sec:deterministic_results}
Numerically, we find two locally optimal solutions, the `suppression' strategy where $\beta$ stays low enough so that ${\bf R_0} <1$ and the `mitigation' strategy where $\beta$ is near $b$ except at the peak of the epidemic, and the epidemic runs its course. This qualitative finding is not sensitive to the choices of parameters. 

\begin{figure}
\centering
\begin{subfigure}[b]{\textwidth}
${\includegraphics[width=\textwidth, trim={0 0 0 2cm}, clip]{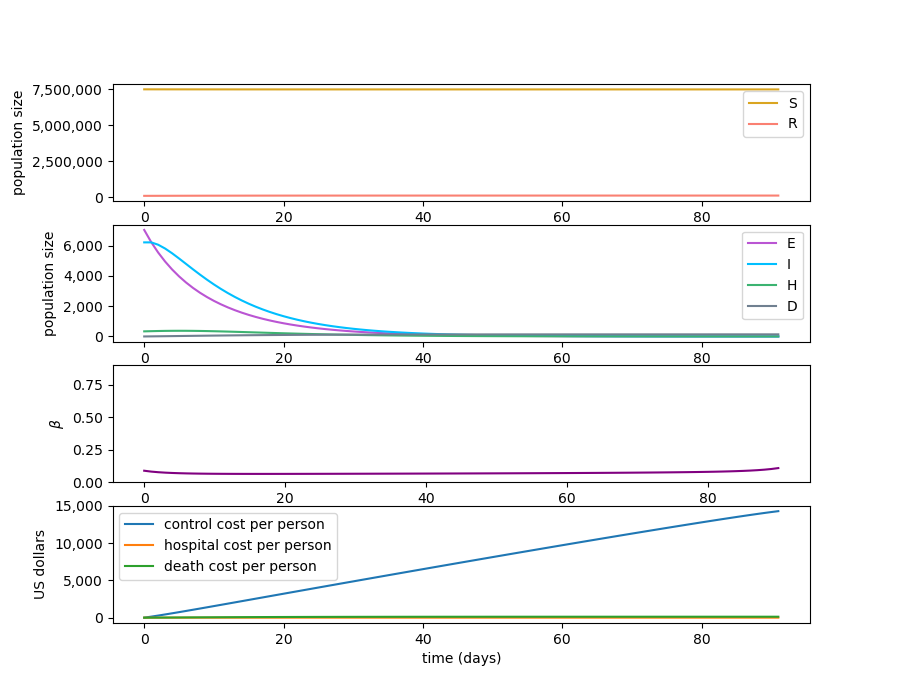}}$
%	\vspace{-1\baselineskip}
        \caption{Supression strategy, with end time 91 days.}
        \label{fig:supression}
    \end{subfigure}
%    	\vspace{-3\baselineskip}
\begin{subfigure}[b]{\textwidth} 
	${\includegraphics[width=\textwidth, trim={0 0 0 2cm}, clip]{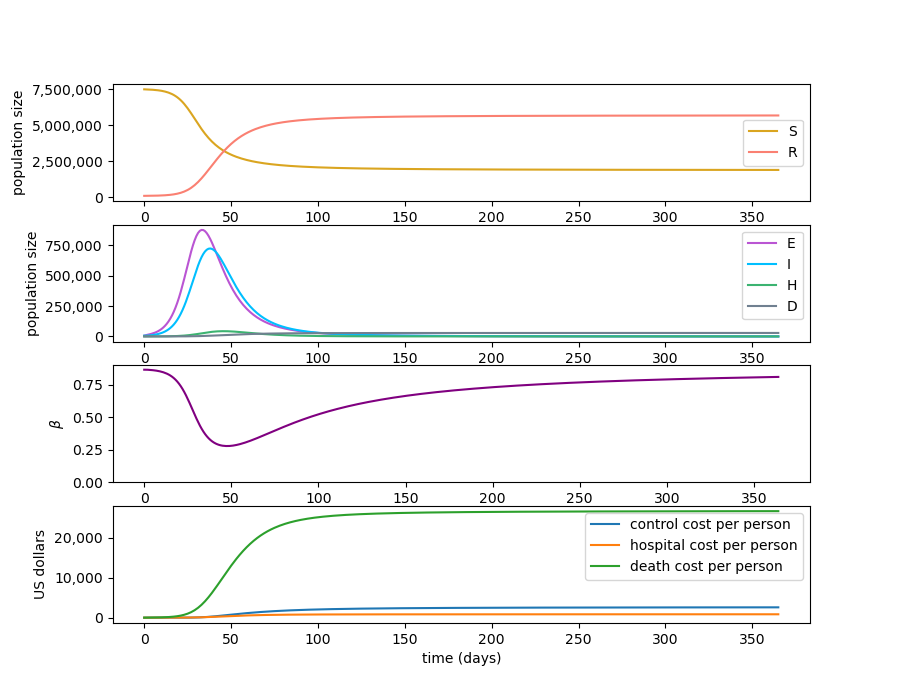}}$
		\caption{Mitigation strategy, with end time 4,061 but only graphed to 365 days.}
		\label{fig:mitigation}
		\end{subfigure}
		\caption{The `suppression' and `mitigation' strategies, top (a) and bottom (b), respectively. Parameters are from Table~\ref{table:parameters} and Table~\ref{table:control_parameters} with initial conditions (\ref{eqn:initial_values}).  See Appendix \ref{apx:deterministic} for details on the numerical approximation.}
\label{fig:optimal_mitigation}
\end{figure}

The two locally optimal solutions, `suppression' strategy and `mitigation' strategy, are plotted in Figure \ref{fig:optimal_mitigation}. The  `suppression' strategy keeps $\beta$ low at a high control cost but low hospitalization and death costs. In contrast, the `mitigation' strategy starts with a high value of $\beta$ (reflecting a relatively open society), and lowers $\beta$ as the disease surges and raises $\beta$ again, resulting in a low control cost and high hospitalization and death costs. The globally optimal solution is the `suppression' strategy, with a cost about a third of the cost of the `mitigation' strategy.  

The cost of the `suppression' strategy is approximately \$15,137 per person, with a total cost of \$115 billion, and is mostly from the {\it control cost}. The cost of the `mitigation' strategy is \$30,226 per person, a total cost of \$229 billion, and mostly from the {\it cost of deaths}.

The end time $T$ for the `mitigation' strategy is much longer than that for the `suppression' strategy. In the `suppression' strategy $T = 91$ and in the `mitigation' strategy $T = 4,061$ (only the first $365$ days are shown in Figure \ref{fig:optimal_mitigation}).  
The long time to elimination for the `mitigation' strategy is due to the end state being very near to $S_T = {N}/{{\bf R_0}}$.
%$S_T = \frac{N}{{\bf R_0}}$.

\begin{figure}
\centering
\begin{subfigure}[b]{0.48\textwidth}
$\includegraphics[width=\textwidth]{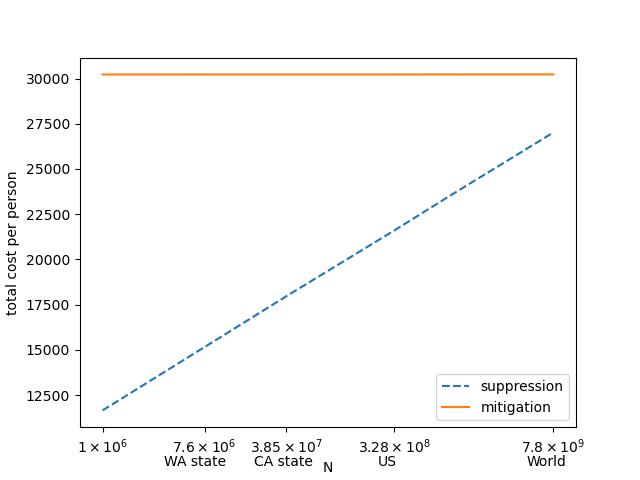}$\ 
\caption{Total cost per person versus
$N$.}
\end{subfigure}
\begin{subfigure}[b]{0.48\textwidth}
$\includegraphics[width=\textwidth]{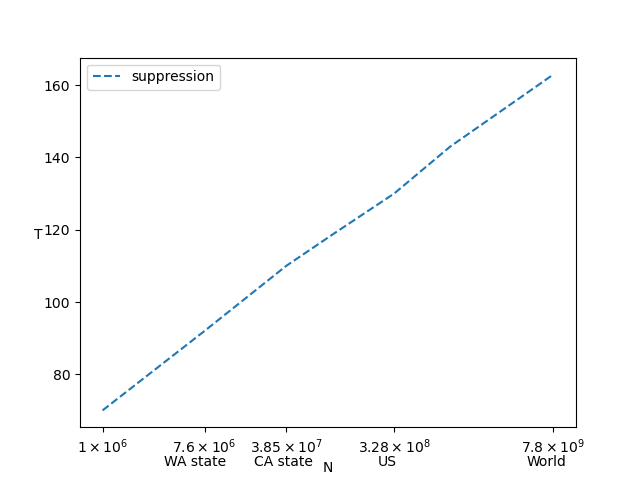}$
\caption{End time $T$ versus $N$.}
\end{subfigure}
\caption{The total cost per person (a)  and  the end time $T$ (b)  for a range of $N$ values from $1$ million to  $7.8$ billion (world population)
 on a log scale.  The mitigation strategy is nearly identical for different $N$ values and its end time is not plotted as it is very large.  The solutions are qualitatively similar across $N$ values.}
\label{fig:stress_test_N}
\end{figure}

The cost of the `suppression' strategy is proportional to $N T$, and $T$ is determined by the end-time threshold. The control variable $\beta$ is kept low so that $I_t, E_t$, and $H_t$ decay exponentially, and the end time is achieved relatively quickly. If the whole model is scaled by $N$, and the end-time threshold remains constant, then the end time $T$ is proportional to $\log N$. 
    The cost of the  `mitigation' strategy is proportional to $N$ but less influenced by the end time, so the  control variable $\beta$ can reach the initial, uncontrolled, infection rate $b$. The total cost for a range of $N$ values, with initial conditions scaled proportional to $N$, is shown in Figure~\ref{fig:stress_test_N}, which demonstrates the above mentioned scaling.
    
The marginal cost of the end-time constraint, as  reflected by the Lagrange multiplier $\sigma$ (see Appendix~\ref{apx:deterministic}), indicates the savings for being allowed to stop while one person is still infected. This marginal cost for the `suppression' strategy  ($\sigma$ is \$2,841,000,000) is much larger than the marginal cost  for the `mitigation' strategy ($\sigma$ is \$2,433,000). The `suppression' strategy terminates much faster than the `mitigation` strategy, so an additional day for the shorter time frame of `suppression' strategy ($T = 91$) has more impact than an additional day for the longer time frame of the `mitigation' strategy ($T = 4,061$).

\begin{figure}
\centering
\begin{subfigure}[b]{0.48\textwidth}
$\includegraphics[width=\textwidth]{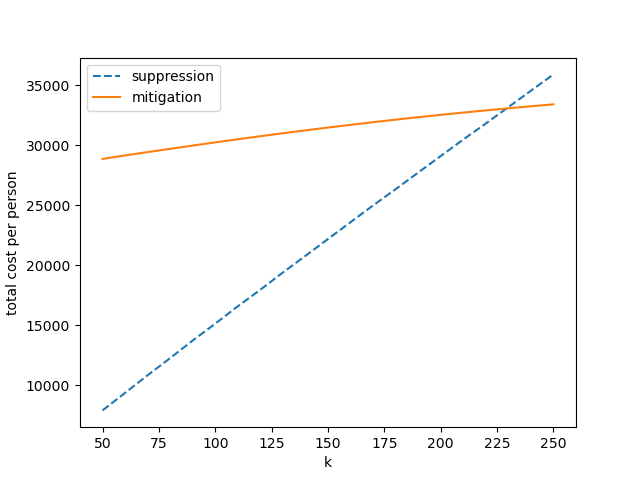}$\ 
\caption{Total cost per person versus $k$.}
\end{subfigure}
\begin{subfigure}[b]{0.48\textwidth}
$\includegraphics[width=\textwidth]{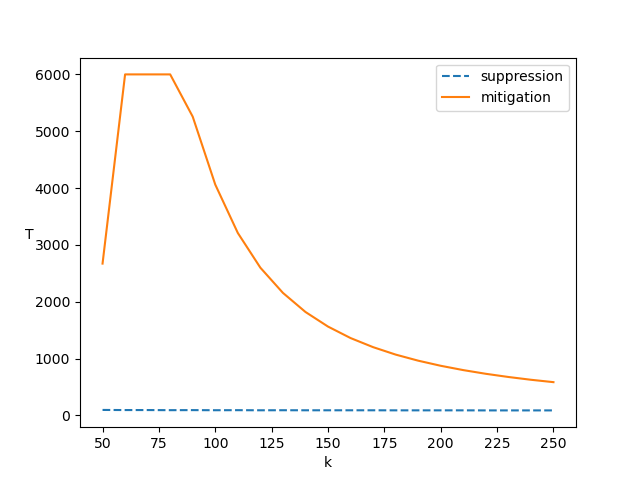}$
\caption{End time $T$ versus $k$.}
\end{subfigure}
\caption{The total cost per person (a) and  the end time $T$ (b) is plotted for a range of $k$ values and the two locally optimal strategies.  The solutions are qualitatively similar across $k$ values, although the end time shows a large variation for the `mitigation' strategy (for the `suppression' strategy the end time simply decreases from 95 to 88 across the $k$ values).  The end time is capped at 6000, which is reached for $k=70$ through $k=90$. The `mitigation' strategy is no longer found as a locally optimal policy at $k=50$.}
\label{fig:stress_test}
\end{figure}

Figure \ref{fig:stress_test} illustrates the impact of a  range of values for the control cost coefficient $k$, ranging from $50$ to $250$, on the cost per person and the end time for both the `suppression' strategy and the `mitigation' strategy.  We find that for a range of $k$ values from 70 to 90, the end time for the `mitigation' strategy effectively diverges to $+\infty$, so numerically, we capped the end time to 6,000 in this range.

\begin{figure}
\centering
$\includegraphics[width=0.5\textwidth]{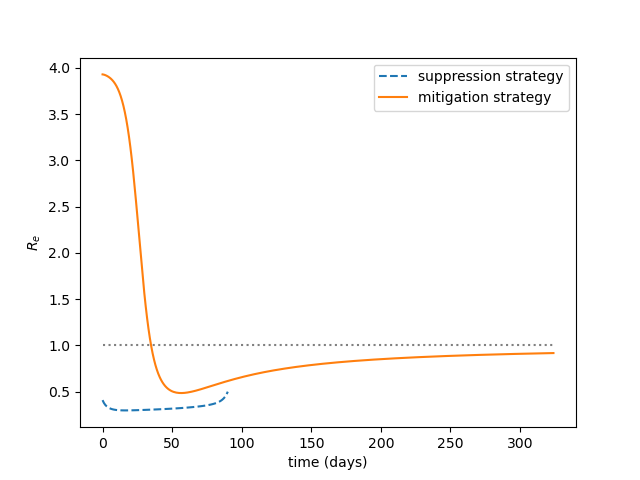}$
\caption{The effective reproduction number in each strategy over time. The  `suppression' strategy ends much sooner than the `mitigation' strategy. {Note that, in the `mitigation' strategy, $R_e$  appears to approach $1$ asymptotically with these parameters.} }
\label{fig:R_e_plot}
\end{figure}

 Figure \ref{fig:R_e_plot} compares the effective reproduction number (as in \eqref{eqn:R_e}) of the two strategies over time.
% For our parameter values, we compute the effective reproduction number as ${\bf R_e}  = \frac{\beta\, \frac{S}{N}}{\lambda_0+\gamma_0+\delta_0}$. 
This gives a useful perspective on the qualitative nature of the two solutions.    We observe that, in the `mitigation strategy', ${\bf R_e}$ appears to approach~$1$ asymptotically with these parameters. 

\begin{figure}
\centering
\begin{subfigure}[b]{0.48\textwidth}
$\includegraphics[width=\textwidth]{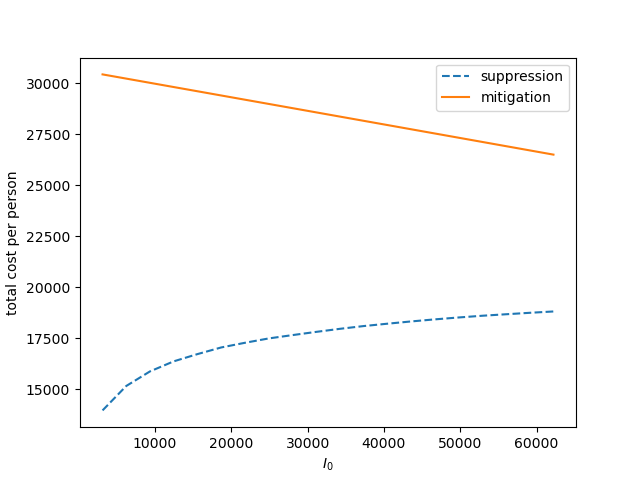}$\ 
\caption{Total cost per person versus initial infectious population $I_0$.}
\end{subfigure}
\begin{subfigure}[b]{0.48\textwidth}
$\includegraphics[width=\textwidth]{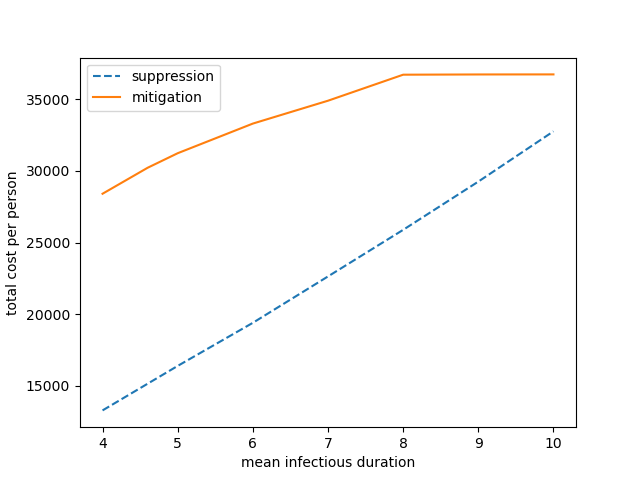}$
\caption{Total cost per person versus mean infectious period.}
\end{subfigure}
\caption{The total cost per person is plotted for a range of $I_0$ values (other initial values scaled accordingly) and a range of the mean infectious period ($1/(\lambda_0+\gamma_0+\delta_0)$ used to determine $\lambda_0$, $\delta_0$, $\gamma_0$).  The solutions are qualitatively similar.}
\label{fig:stress_test_I_lambda}
\end{figure}

Figure \ref{fig:stress_test_I_lambda} provides a final `stress test' of some of our parameters, plotting the optimal cost for both strategies for a range of initial conditions (where $E_0$, $H_0$, and $R_0$ are scaled proportional to $I_0$) and for a range of values for the mean infectious duration  (corresponding to $(\lambda_0+\gamma_0+\delta_0)^{-1}$ with $\lambda_0$, $\gamma_0$, $\delta_0$ scaled accordingly).  We find the solutions remain qualitatively similar and the preference of the `suppression' strategy as the global optimum remains for a wide range of these parameters.

\section{Vaccinations}
\label{Sec:vaccinations}

We now consider the effect of vaccination on optimal control of the COVID-19 pandemic.  We introduce an additional transition from the susceptible state, $S$, to the recovered state, $R$, representing an individual gaining an immunity to COVID-19 by vaccination. See illustration in Figure~\ref{fig:ratediagram}.  We assume that the vaccine is distributed at a constant rate, determined by the parameter $o$ for the ratio of population vaccinated per day.  Thus the equation for the susceptible population, when $S_t>0$, becomes
$$
  \frac{dS_t}{dt} = -\, \beta_t\, S_t\, \frac{I_t}{N}\, -\, o\, N.
$$
Clearly, this results in $S_t$ reaching $0$ in a finite time (before $t=o^{-1}$) , after which we assume that $S_t$ remains at $0$.  The only other modification of the model to include vaccinations is the addition of the $o\, N$ term in the recovered population, 
$$
    \frac{dR_t}{dt} = \gamma_0\, I_t\, +\gamma_1\, H_t+\, o\, N.
$$
\subsection{Vaccination Parameters}

The rate of vaccination, $o\, N$, depends on the parameter $o$, where $o^{-1}$ represents the number of days to vaccinate the entire population.
We consider the parameter values ranging from $o=1/250$, reflecting if vaccination begins on January 1, 2021 then the entire population is vaccinated by September 8, 2021,  to $o=1/500$, where if vaccination begins on January 1, 2021 then the entire population is vaccinated by May 16, 2022.

We also consider a second set of initial conditions, which are chosen to represent COVID-19 in the U.S.~as of January 1, 2021, with  $N=328.2$ million, and
\begin{align}\label{eqn:initial_values_US}
	\left(\begin{array}{c}
		S_0\\
		E_0\\
		I_0\\
		H_0\\
		R_0\\
		D_0
	\end{array}\right) = 
	\left(\begin{array}{c}
		328,200,000 - E_0 -I_0 - H_0 - R_0\\
		190,728 * 4.6 / \alpha\\
		190,728 * 4.6 / (\lambda_0+\gamma_0+\delta_0)\\
		125,047 * 1.9\\
		20,112,544 * 4.6 - I_0 - E_0 - H_0\\
		0
	\end{array}\right)
	=\left(\begin{array}{c}
		235,682,298\\
		4,569,525\\
		4,035,804\\
		237,589\\
		83,674,784\\
		0
	\end{array}\right).
\end{align}
We  used 190,728 for the 7-day rolling average of daily confirmed cases and 20,112,544 total cases \cite{CDC-tracker} on Jan 1, 2020.  We  used 125,047  current hospitalizations as of Jan 1, 2021 \cite{covid-tracking-project}. We again used the estimate of 4.6 times more cases than confirmed cases, and 1.9 times more hospitalized cases than reported \cite{CDC-burden}.

\subsection{Vaccination Results}

We find a new optimal strategy, which we call `delay-mitigation,' that delays the epidemic by maintaining disease levels at a near constant state, until the majority of the population has been vaccinated and the infection dies out in the population due to immunity.  This new strategy takes the place of the `mitigation' strategy, because the end time for the `mitigation' strategy greatly exceeds the time to vaccinate the entire population, $o^{-1}$. 

\begin{figure}
\centering
\begin{subfigure}[b]{\textwidth}
$\includegraphics[width=\textwidth, trim={0 0 0 2cm}, clip]{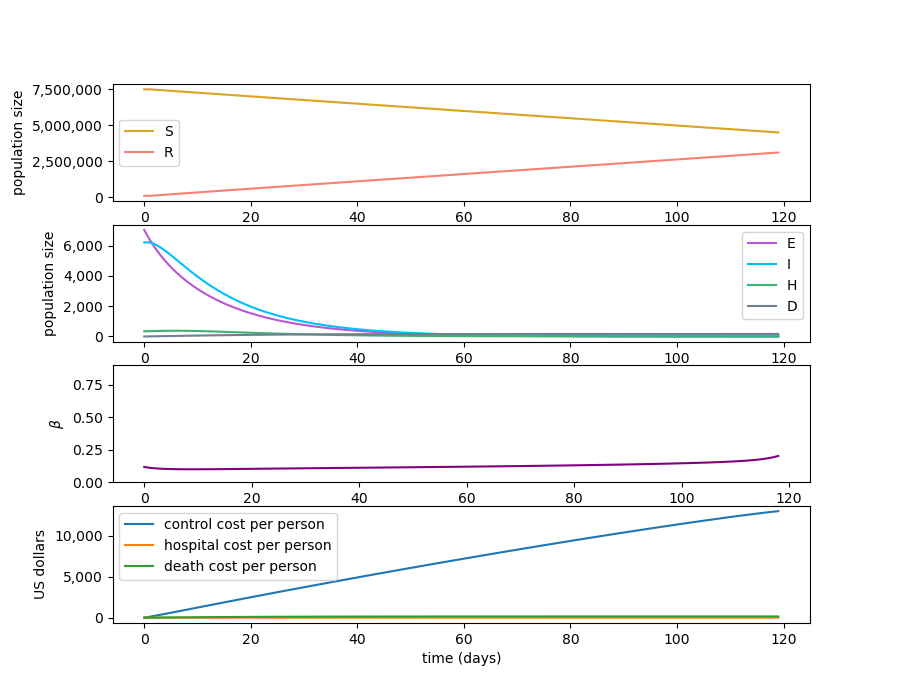}$
        \caption{Supression strategy, with vaccine.}
        \label{fig:supression_vaccine}
    \end{subfigure}
\begin{subfigure}[b]{\textwidth} 
	$\includegraphics[width=\textwidth, trim={0 0 0 2cm}, clip]{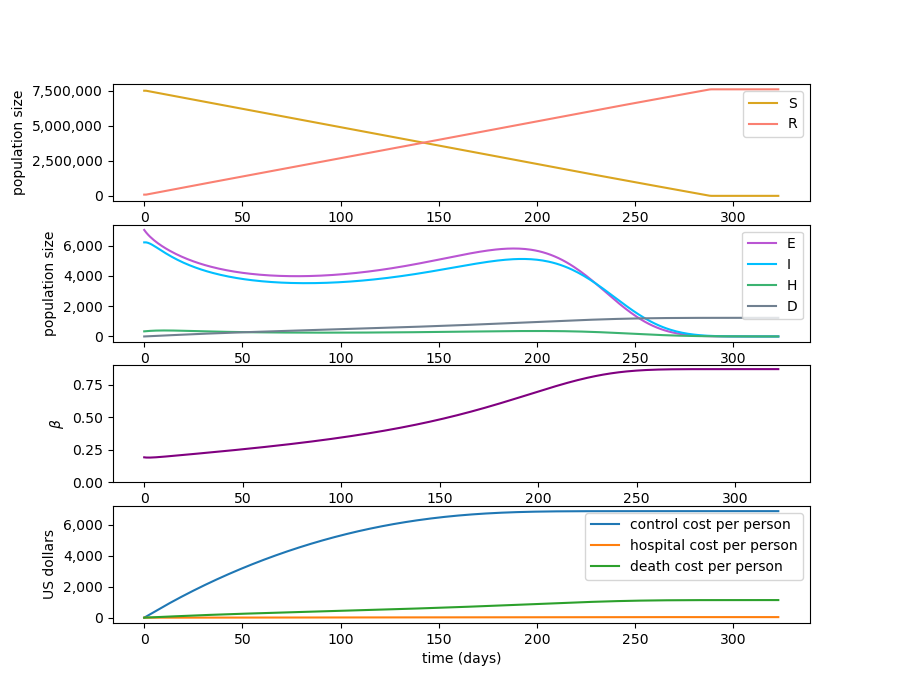}$
		\caption{Delay-mitigation strategy, with vaccine.}
		\label{fig:mitigation_vaccine}
		\end{subfigure}
		\caption{The `suppression' strategy (a) and `delay-mitigation' strategy (b)  with vaccinations. Parameters are from Table~\ref{table:parameters} and Table~\ref{table:control_parameters} with $o=1/300$ and initial conditions (\ref{eqn:initial_values}) reflecting Washington State.}
\label{fig:delay_mitigation_vaccine}
\end{figure}

For the initial conditions (\ref{eqn:initial_values}) reflecting Washington State, and values of $o^{-1}\geq 300$ (i.e., $o\leq 1/300$), we still find a `suppression' strategy as a local optima that is nearly identical to the `suppression' strategy found in Section \ref{sec:deterministic_results}. However, the `delay-mitigation' strategy is the global optima. When $o=1/250$, the `suppression' strategy is not locally optimal, because it is better to vaccinate more of the population by extending the end time.   With $o=1/300$ and initial conditions (\ref{eqn:initial_values}), the `suppression' and `delay-mitigation' strategies are shown in Figure \ref{fig:delay_mitigation_vaccine}.  

\begin{figure}
\centering
$\includegraphics[width=\textwidth]{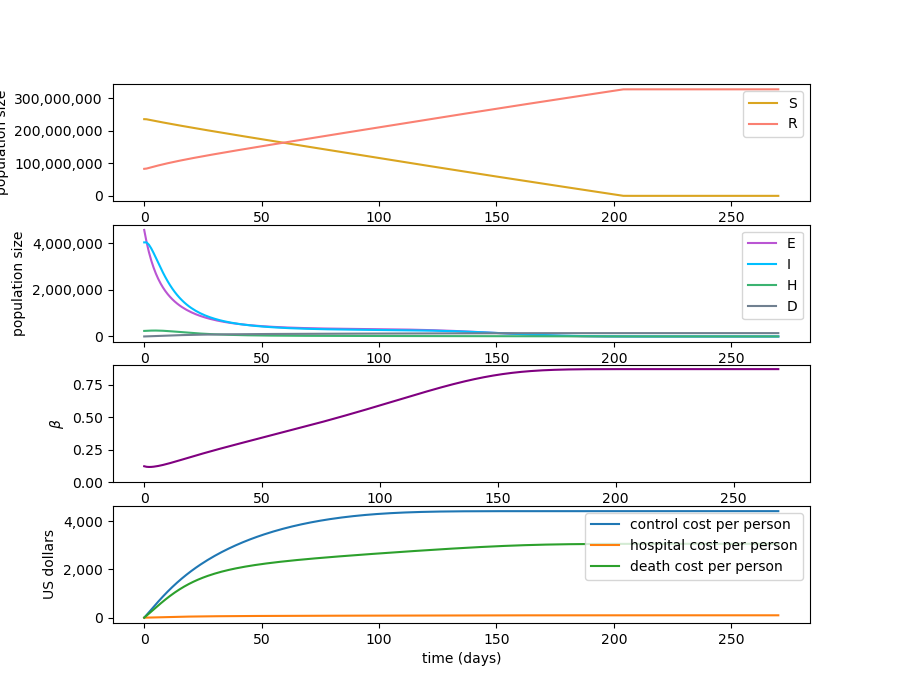}$
\caption{The `delay-mitigation' strategy is the only optimal strategy when we choose initial conditions \eqref{eqn:initial_values_US} reflecting the US population at the beginning of 2021 ($N=328,200,000$), with $o=1/300$.}
\label{fig:optimal_mitigation_vaccine_US}
\end{figure}

For the initial conditions (\ref{eqn:initial_values_US}) reflecting the U.S., we never find the `suppression' strategy to be locally optimal, and the only local optimal strategy is the `delay-mitigation', which is shown in Figure \ref{fig:optimal_mitigation_vaccine_US}.

% This `delay-mitigation' strategy appears to be the global minimum for $o\gtrsim 1/525$. The total cost per person of both strategies are plotted for a range of $o$ values in Figure \ref{fig:o_stress_test}.

% With vaccination, the `mitigation' strategy no longer exists as the end-time cannot be much longer then the time to vaccinate the entire population, $o^{-1}$.

% For values of $o\lesssim1/275$ we find a `suppression' strategy nearly identical to the `suppression' strategy found in Section \ref{sec:deterministic_results}.  When $o=1/250$, this solution appears to not be locally optimal.  

% Rather than a `mitigation' strategy, we find a locally optimal strategy that instead delays the epidemic, until the majority of the population has been vaccinated and the infection dies out in the population due to immunity.  This `delay-mitigation' strategy appears to be the global minimum for $o\gtrsim 1/525$. The total cost per person of both strategies are plotted for a range of $o$ values in Figure \ref{fig:o_stress_test}.

%\left(\begin{array}{c}S_0 \\ E_0 \\I_0\\ H_0\\ R_0\\ D_0\end{array}\right)=\left(\begin{array}{c}245,000,000 \\ 1,322,711 \\1,170,000\\ 82,900\\ 80,124,389\\ 400,000\end{array}\right).
%\end{align}
%{\red Hospitalization is maybe a bit low... but this is consistent with our parameters.  Maybe deaths should be 500,000 to include excess deaths?}

\begin{figure}
\centering
\begin{subfigure}[b]{\textwidth}
      \begin{subfigure}[b]{0.48\textwidth}
      $\includegraphics[width=\textwidth]{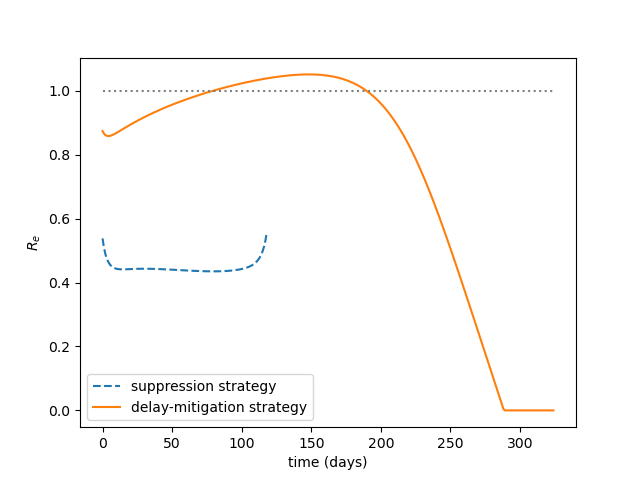}$\ 
      \caption{Effective reproduction number in each strategy with initial conditions of (\ref{eqn:initial_values}). }
\end{subfigure}
\begin{subfigure}[b]{0.48\textwidth}
$\includegraphics[width=\textwidth]{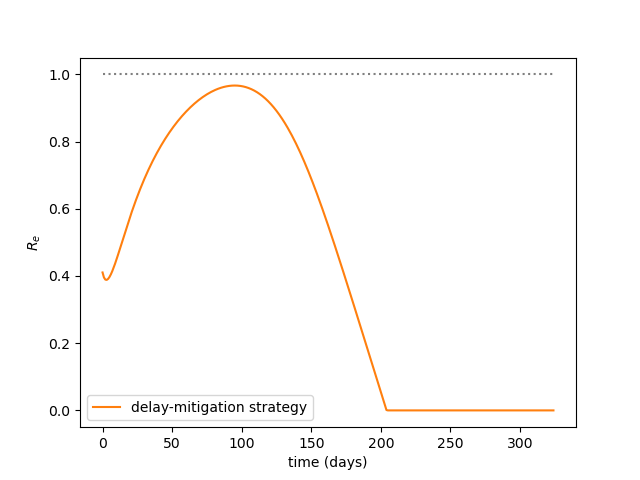}$
\caption{Effective reproduction number  with initial conditions of (\ref{eqn:initial_values_US}).}
\end{subfigure}
\label{fig:R_e_plot_vaccine}
\end{subfigure}
%NEXT TWO FIGURES
\begin{subfigure}[b]{\textwidth}
\centering
\begin{subfigure}[b]{0.48\textwidth}
$\includegraphics[width=\textwidth]{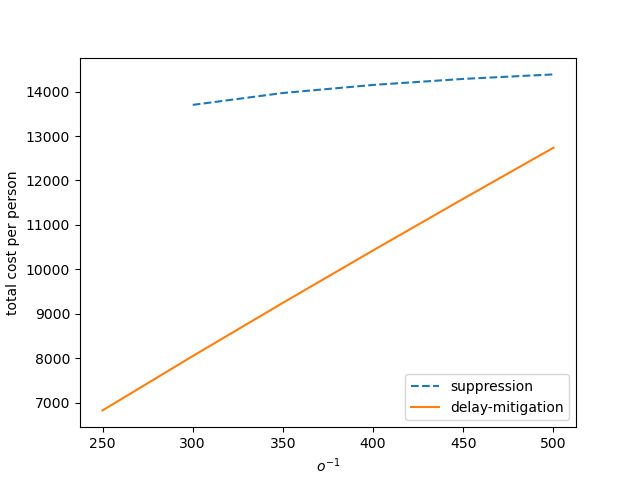}$ \
\caption{Total cost per person versus $o^{-1}$ with initial conditions  (\ref{eqn:initial_values}). }
\end{subfigure}
\begin{subfigure}[b]{0.48\textwidth}
$\includegraphics[width=\textwidth]{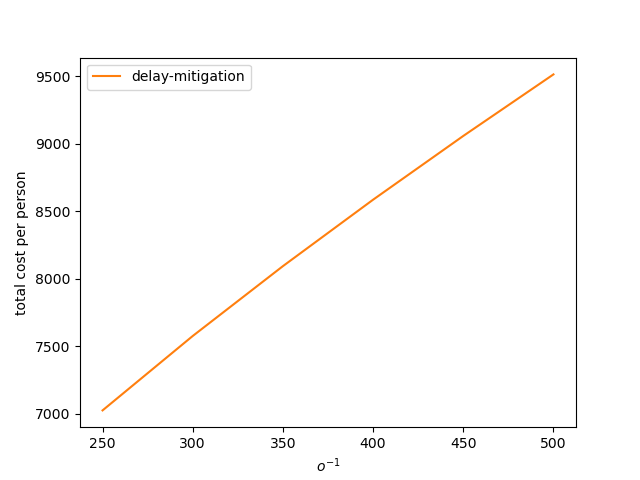}$
\caption{Total cost per person versus $o^{-1}$ with initial conditions  (\ref{eqn:initial_values_US}).}
\end{subfigure}
\label{fig:o_stress_test}
\end{subfigure}
%\caption{The effective reproduction number  shown with initial conditions of (\ref{eqn:initial_values}) in (a) and of (\ref{eqn:initial_values_US}) in (b).  For initial conditions of (\ref{eqn:initial_values}) corresponding to Washington State (a), both strategies are found as local optima, whereas, for  initial conditions of (\ref{eqn:initial_values_US}) corresponding to the US in early 2021 (b), we do not find a `suppression' strategy as a local optima and instead the `delay-mitigation' strategy incorporates strong control of the infection rate in early stages. 
%The total cost per person for ranges of time to vaccinate the entire population $o^{-1}$, ranging from 250 to 500  ($o$ values ranging from $1/250$ to $1/500$).  The initial conditions of (\ref{eqn:initial_values}) reflecting Washington state are in (c) and of (\ref{eqn:initial_values_US}) reflecting the US  in (d). }
\caption{The effective reproduction number  shown with initial conditions (\ref{eqn:initial_values}) in (a) and  (\ref{eqn:initial_values_US}) in (b). The total cost per person for ranges of time to vaccinate the entire population $o^{-1}$, ranging from 250 to 500  ($o$ values ranging from $1/250$ to $1/500$), with initial conditions   (\ref{eqn:initial_values}) reflecting Washington state in (c), and with initial conditions (\ref{eqn:initial_values_US}) reflecting the US  in (d).}
\label{fig:Four-part_delay_mitigation_vaccine}
\end{figure}

For the parameter value $o=1/300$, we find:
\begin{itemize}
    \item With initial conditions (\ref{eqn:initial_values}), the cost of the `suppression' strategy is \$13,701 per person, down only slightly with the introduction of the vaccine, and the cost of the `delay-mitigation' strategy is \$8,041 per person, with the cost mostly from the control cost (compared with the `mitigation' strategy, which had cost mostly due to deaths).  However, the end time for the the `suppression' strategy is less than half that of the `delay-mitigation' strategy.  In the  `suppression' strategy, the disease ends with nearly half the population in the Susceptible compartment, whereas in the  `delay-mitigation' strategy,  the entire population has left the Susceptible compartment by the end time.   In the  `suppression' strategy, about half the population is vaccinated by the time the disease is under control. The dynamics, control, and cost are plotted for the `suppression' strategy and the `delay-mitigation' strategy in Figure~\ref{fig:delay_mitigation_vaccine} with initial conditions (\ref{eqn:initial_values}).
    
    \item With initial conditions (\ref{eqn:initial_values_US}), only the `delay-mitigation' strategy is found to be optimal, with total cost \$7,556 per person.  We note the cost per person is slightly lower than that using initial conditions (\ref{eqn:initial_values}), even though the disease levels are higher, because a larger proportion of the population has already recovered from the disease.   The dynamics, control, and cost are plotted for the single optimal, `delay-mitigation' strategy with initial conditions (\ref{eqn:initial_values_US}) in Figure~\ref{fig:optimal_mitigation_vaccine_US}. 
    
    \item The end time $T$ for the `delay-mitigation' strategy is $T=323$ and $T=270$, for initial conditions (\ref{eqn:initial_values}) and (\ref{eqn:initial_values_US}), respectively, finishing after the entire population has been vaccinated and the last infected person has recovered. In the `suppression' strategy,  the end time is shorter, $T = 119$, which is similar but slightly longer than the end time without a vaccine.

    \item We observe a new shape for the control policy in the `delay-mitigation' strategy, as is evident in the shape of  the effective reproduction number, which is maintained near 1, until the majority of the population has been vaccinated. A plot of ${\bf R_e}$ is in Figure~\ref{fig:Four-part_delay_mitigation_vaccine}(a) and (b), for different initial conditions.  As the remainder of the population receives a vaccine we observe the effective reproduction number falls to $0$.   When we use the initial conditions (\ref{eqn:initial_values_US}), we see that ${\bf R_e}$ remains below $1$ due to the higher disease level in the initial population.

    \item The total cost per person of the `suppression' strategy only varies slightly with the vaccine rate $o$, while the total cost per person of the `delay-mitigation' strategy appears to depend linearly on $o^{-1}$ as observed in Figure~\ref{fig:Four-part_delay_mitigation_vaccine}(c) and~(d).
\end{itemize}

% \begin{figure}
% \centering
% $\includegraphics[width=\textwidth]{Figures/s_I_o300}$
% \caption{The `suppression' strategy solution with vaccinations. Parameters are from Table~\ref{table:parameters} and Table~\ref{table:control_parameters} with $o=1/300$. We see the disease is eliminated before the vaccine has much impact.}
% \label{fig:optimal_suppression_vaccine}
% \end{figure}

%\begin{figure}
%\centering
%$\includegraphics[width=\textwidth]{Figures/s_o300}$
%\label{fig:suppression_vaccine}
% \end{figure}
% \begin{figure}
% \centering
%$\includegraphics[width=\textwidth]{Figures/m_o300}$
%\caption{The `suppression' strategy and `delay-mitigation' strategy solution with vaccinations. Parameters are from Table~\ref{table:parameters} and Table~\ref{table:control_parameters} with $o=1/300$ and initial conditions (\ref{eqn:initial_values}) reflecting Washington State.}
%\label{fig:delay_mitigation_vaccine}
%\end{figure}
%%%%%

%%TESTING - PUT in ONE 4-part Figure%%

\section{Conclusion}

We have analyzed the optimal control of an idealized epidemic model of COVID-19.  We found two locally optimal strategies, `suppression' and `mitigation', which correspond to qualitatively distinct approaches to combat the epidemic, and the `suppression' strategy is the global optimum.  By considering a wide range of parameters we find the solutions to be fairly robust.  The optimal control strategies provide insight into how the strength or relaxation of control impacts the population and associated costs over time, as well as the impact of the end time.   In the `suppression' strategy, the susceptible population remains high, while in the  `mitigation' strategy, nearly the entire population contracts the disease.  We also see that the end time for the `suppression' strategy is much shorter than the end time for the `mitigation' strategy.

When accounting for vaccinations, we find a globally optimal `delay-mitigation' strategy that delays the spread of the disease until the majority of the population has received a vaccine. The `susceptible' strategy is not globally optimal when vaccination is present, and only appears as a local optimum under certain parameter values. We vary the time it takes to vaccinate the entire population and initial conditions and observe the associated change in the effective reproduction number and the total cost per person.

There are many additional features that we have not attempted to model.  One feature is that the control cost will likely depend on the number of infected individuals.  In particular, targeted contact tracing and quarantine may serve to reduce the infection when the number of infections is small with less cost than overarching social distancing policies.  Similarly, targeted vaccination can also effectively reduce the infection rate as well as remove individuals from the susceptible population.

We have not accounted for the possibility of mutations of the virus that could change the parameters and reduce the effectiveness of immunity in the recovered population.  We consider a stochastic model (in the appendix) as a starting point to analyze the uncertainty introduced by mutations, which would be an important area for further research.

Another feature is the network dependence of epidemic spread, either through social networks or geographic distance.  This is an active area of research with many different existing approaches.

A final feature to mention is the role of information.  Our idealized model has assumed perfect information about the state of the disease, which is not realistic. 
Gaining accurate information about the  parameters and progression of the disease is needed for implementing optimal epidemic control.

\section*{Acknowledgments} This work has been funded in part by the National Science Foundation grant CMMI-1935403.

   \appendix

   \section{Equilibria Analysis}\label{apx:equilibria}

We investigate the Jacobian matrix and its eigenvalues to parametrically characterize the system stability near equilibrium. {If all eigenvalues are less than or equal to zero, the system is stable.} {The zero eigenvalues represent the degrees of freedom of equilibria, which are simply $S$, $R$, and $D$, and do not complicate the analysis.}

The Jacobian matrix associated with the dynamics in \eqref{eqn:mf_de} is 
$$
    \left(\begin{array}{cccccc}  
        -\frac{\beta\, I}{N} & 0& -\beta\, \frac{S}{N}  & 0 & 0 & 0\\
        \frac{\beta\, I}{N} & -\alpha & \beta\, \frac{S}{N} &  0 & 0 & 0\\
        0 &  \alpha & -(\lambda_0+\gamma_0+\delta_0) & 0 & 0 & 0 \\
        0 & 0 & \lambda_0 &  -\gamma_1 -\delta_1 & 0 & 0\\
        0 & 0 & \gamma_0 & \gamma_1 & 0 & 0\\
        0 & 0 & \delta_0 & \delta_1 & 0 & 0 
    \end{array}\right).
$$
At an equilibrium point, the first column becomes zero because $I=0$, and we can compute the eigenvalues  $\epsilon$ by finding the roots of the characteristic polynomial, solving 
$$
   0=det \left(\begin{array}{cccccc}
        -\epsilon & 0& -\beta\, \frac{S}{N}  & 0 & 0 & 0\\
        0 & -\alpha -\epsilon & \beta\, \frac{S}{N} &  0 & 0 & 0\\
        0 &  \alpha & -(\lambda_0+\gamma_0+\delta_0) - \epsilon & 0 & 0 & 0 \\
        0 & 0 & \lambda_0 &  -\gamma_1 -\delta_1  - \epsilon& 0 & 0\\
        0 & 0 & \gamma_0 & \gamma_1 &  -\epsilon & 0\\
        0 & 0 & \delta_0 & \delta_1 & 0 & -\epsilon 
    \end{array}\right),
$$
which yields
$$
0=    \epsilon^3(\gamma_1+\delta_1+\epsilon)\Big((\alpha+\epsilon)(\lambda_0+\gamma_0+\delta_0+\epsilon) - \alpha\, \beta\,  \frac{S}{N}\Big).
$$

There is a zero eigenvalue with multiplicity three, $-\gamma_1-\delta_1$ is an eigenvalue,  and the other two eigenvalues solve
$$
 0=   \epsilon^2+\big(\alpha+[\lambda_0+\gamma_0+\delta_0]\big)\epsilon+\alpha \, [\lambda_0+\gamma_0+\delta_0]-\alpha\, \beta\, \frac{S}{N}
$$
so
$$
    \epsilon = \frac{-\big(\alpha+[\lambda_0+\delta_0+\gamma_0]\big)\pm \sqrt{\big(\alpha+[\lambda_0+\delta_0+\gamma_0]\big)^2-4\, \alpha\, \big( [\lambda_0+\delta_0+\gamma_0]-\beta\, \frac{S}{N}\big)}}{2}.
$$
Both of these eigenvalues are negative if and only if
$$
    \lambda_0+\delta_0+\gamma_0>\beta\, \frac{S}{N}.
$$
{This eigenvalue analysis provides a condition on parameters that ensure  system stability.}

\begin{figure}
\centering
$
\includegraphics[width=0.7\textwidth]{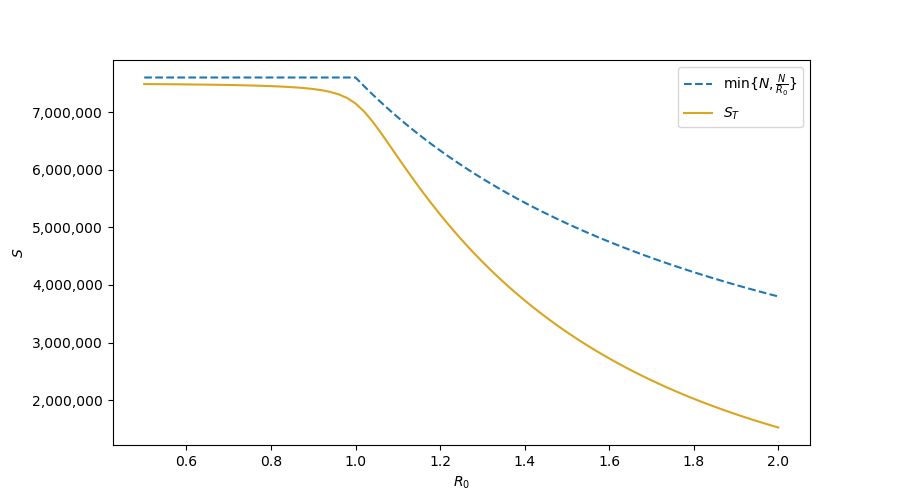}
$
\caption{The final state $S_T$ as function of ${\bf R_0}$, a theoretical upper bound in the dashed blue line and simulations  in yellow. Parameters are as in Table~\ref{table:parameters} with $\beta$ ranging as a function of ${\bf R_0}$, $\beta = (\lambda_0 + \gamma_0+\delta_0)\, {\bf R_0}$. The end time for the simulations is taken to be $T = $ 6,000 days.
\label{fig:R_0}}
\end{figure}

We plot the equilibrium value of $S$ at a large time $T$ for a range of ${\bf R_0}$ values in Figure~\ref{fig:R_0}.
% as a result of simulation and the theoretical bound of Theorem~\ref{thm:equilibrium}.  
The ${\bf R_0}$ value of COVID-19 without interventions has been estimated to be greater than 2, which would result in the majority of the population becoming infected.  A clear qualitative feature shown in Figure~\ref{fig:R_0} is that $S_T$ is near $N$ when ${\bf R_0}\leq1$ and the final state of $S_T$ always lies below an upper bound of $\frac{N}{{\bf R_0}}$.  We find in Section~\ref{sec:deterministic_results} that the optimal `mitigation' strategy achieves the end result of $S_T \approx \frac{N}{{\bf R_0}}$.     

 \section{Pontryagin Maximum Principle}\label{apx:pontryagin}

We recall the dynamics (\ref{eqn:mf_de}) and define $(\Sigma)_{t\in [0,T]} = (S_t, E_t, I_t, H_t, R_t, D_t)$ and let $f(\Sigma,\beta)\in \R^d$ denote the righthand side, i.e., 
\begin{align}\label{eqn:sigma_dynamics}
    \frac{d\Sigma_t}{dt}=f(\Sigma_t,\beta_t).
\end{align}

We let $(P_t)_{t\in [0,T]}=(P^S_t,P^E_t,P^I_t,P^H_t,P^R_t,P^D_t)_{t\in [0,T]}$ be the costate and define the Hamiltonian
$$
    \mathcal{H}(\Sigma,P,\beta) =  f(\Sigma,\beta)\cdot P-L(\beta)-F(H).
$$
Then we let $(P_t)_{t\in [0,T]}$ solve the costate equation, applying a differential operator to the Hamiltonian,
\begin{align}\label{eqn:p_dynamics}
    \frac{dP_t}{dt} = - D_\Sigma \mathcal{H}(\Sigma_t,P_t,\beta_t),
\end{align}
with
$$
    \left(\begin{array}{c} P^S_T \\ P_T^E + \sigma \\ P^I_T + \sigma \\ P^H_T + \sigma \\ P^R_T \\ P^D_T\end{array}\right) = -D_\Sigma G(\Sigma_T).
$$
Here $\sigma \geq 0$ is a Lagrange multiplier for the target constraint.  The free end time yields the additional transversality condition
\begin{align}\label{eqn:transversality}
    \sup_{\beta}\mathcal{H}(\Sigma_T, P_T, \beta) = 0.
\end{align}

The Pontryagin Maximum Principle states:
\begin{theorem} \label{thm:pontryagin}
If $T$ and $(\beta_t)_{t\in [0,T]}$ are optimal and $(\Sigma_t)_{t\in [0,T]}$ solves {\normalfont(\ref{eqn:mf_de})} with {\normalfont(\ref{eqn:extinction})} satisfied, then there is $\sigma\geq 0$ and $(P_t)_{t\in [0,T]}$ that solves {\normalfont(\ref{eqn:p_dynamics})} with terminal conditions, such that $\sigma\, (I_T + E_T + H_T - e^{-1}) =0$, the transversality condition {\normalfont (\ref{eqn:transversality})} is satisfied, and for almost every $t\in [0,T]$,
$$
    \beta_t\in {\rm argmax}\{\mathcal{H}(\Sigma_t,P_t,\cdot)\}.
$$
\end{theorem}

We note that  $\beta_t\in {\rm argmax}\{\mathcal{H}(\Sigma_t,P_t,\cdot)\}$ provides optimal $\beta_t$ given the state $\Sigma_t$, and co-state $P_t$. The Pontryagin Maximum Principle provides necessary, but not sufficient conditions for optimality.   We numerically discover two local optima that satisfy these conditions, that we call `suppression' and `mitigation' strategies.

We let $J_\sigma$ denote the augmented cost 
$$
    J_\sigma\big[(\beta_s)_{s\in [0,T]}, T\big] = J\big[(\beta_s)_{s\in [0,T]}, T\big] + \sigma\, \big(I_T + E_T + H_T - e^{-1}\big),
$$
and note that for $\sigma$ of Theorem \ref{thm:pontryagin}, $T$ and $(\beta_t)_{t\in [0,T]}$ minimize $J_\sigma$ over policies unconstrained by the end-time threshold (\ref{eqn:extinction}).
For any smooth $(\beta_t)_{t\in [0,T]}$, and solutions $(\Sigma_t)_{t\in [0,T]}$ and $(P_t)_{t\in [0,T]}$, we can calculate the functional derivative
\begin{align}\label{eqn:J_beta_derivative}
    D_{(\beta_s)_{s\in [0,T]}} J_\sigma\big[(\beta_s)_{s\in [0,T]}, T\big](t) = -D_\beta \mathcal{H}(\Sigma_t,P_t,\beta_t),
\end{align}
and
\begin{align}\label{eqn:J_T_derivative}
    D_TJ_\sigma\big[(\beta_s)_{s\in [0,T]}, T\big] = - \mathcal{H}(\Sigma_T,P_T, \beta_T).
\end{align}

Note that even if $L$, $F$, and $G$ are convex, our problem is not  convex due to the nonlinear dynamics.  We use a discretized version of the functional gradients to search for local optima, and compare the overall cost of the strategies to determine the global optimum. 

   \section{Deterministic Control Numerics}\label{apx:deterministic}
    
    We discretize time in fixed increments of $\Delta t = 1$.  Dynamics are approximated by a first order Euler scheme, which allows for an easy and exact calculation of the discrete gradients corresponding to (\ref{eqn:J_beta_derivative}) and (\ref{eqn:J_T_derivative}).  The target constraint is relaxed by adding the quadratic penalty function to the cost:
    $$
        \frac{N}{2\mu} \big(max\{0,E_T + I_T + H_T - e^{-1}\}\big)^2
    $$
    and the Lagrange multiplier is retrieved simply as $\sigma = \frac{N}{\mu}\big(max\{0,E_T + I_T + H_T - e^{-1}\}\big)$. We use $\mu = 0.01$.

    We then calculate the gradient of the cost with respect to the control variable by back propagation (equivalent to (\ref{eqn:J_beta_derivative}) and a discretization of costate equations) and employ momentum gradient descent (momentum factor is $0.9$).  The gradient step is chosen to be around $10^{-4}/ N\sim 10^{-5}/N$ for different parameter values. Larger gradient steps lead to instabilities and slow convergence, especially when there is a long time horizon.

    To find the optimal end time we implement the simple algorithm that if the Hamiltonian at the current end time $T$ is positive, we increase the end time by one increment, and otherwise, if the Hamiltonian at $T - \Delta t$ is negative we decrease the end time to $T - \Delta t$.  We find this approach works numerically. 
    
  \section{Stochastic Epidemic Control}{\label{apx:stochastic}

\subsection{Stochastic Considerations}

{We  presented a deterministic model, with continuously varying state variables, as an approximation to a stochastic model, which would consider the transitions between states as random jumps corresponding to individual infections/recoveries/etc.  When the population is large, the transition of a single individual is small relative to the number of individuals in the state.  However, there are features of the stochastic model that are not captured in the deterministic model.  One of these features is the distribution of the end time.}   {As a proof of concept, we now consider the stochastic version of the deterministic models in the paper.}  This illustrates the possibility of including other uncertain aspects such as the introduction of a new strain of COVID-19 with uncertain characteristics.

The error of a deterministic approximation to a stochastic model at a given time is proportional to the standard deviation of the macroscopic state variables, which is of order $\sqrt{N}$, so that the relative error, of order $\frac{1}{\sqrt{N}}$, becomes small when $N$ is large.      
    The time to reach equilibrium is finite in the stochastic model, which we have incorporated into the deterministic model by ending when the infected population reaches a given threshold, $e^{-1}$. If ${\bf R_0}<1$, the end time is of order $\log(N)$, and the error of the deterministic approximation is again proportional to the standard deviation, which is of order $\sqrt{\log(N)}$.  At least for the `suppression' type strategy, the cost is of order $N\, T \sim N\, \log(N)$ and the relative error grows if we scale by $N$.  

    A second interesting note is that when approaching this problem from a dynamic programming perspective, the discretization of the deterministic model naturally leads to a stochastic interpretation by addition of numerical viscosity.  We  proceed to take this dynamic programming perspective and properly account for the fluctuations of the stochastic SEIHRD model.

    An alternative approach to account for the stochastic fluctuations, which we do not consider here, is to approximate the problem near the deterministic solution as a linear quadratic Gaussian (LQG) stochastic optimal control problem.  Since the dynamics of the SEIHRD model are not linear, this approximation must be done carefully. The LQG problem can then be solved as a system of Ricatti differential equations.  It is not clear if the LQG approach can account for the fluctuations of the end time since it is usually done with a fixed end time, but it has proven very effective in many applications and can also account for noisy and incomplete observations of the state variables.

\subsection{Stochastic Simplified Problem}\label{sec:stochastic_simplified}
We  now use the dynamic programming approach to better understand the stochastic nature of the epidemic model.  This approach also has the feature that it finds the globally optimal strategy, where as the approach of Section \ref{sec:det} centered around finding local optima.

The dynamic programming approach is not practical to solve with six state variables, as it would require a computational complexity on the order of $N^6$.  Instead we reduce to a simplified three state model $\tilde{S}, \tilde{I}, \tilde{R}$, as depicted in Figure~\ref{fig:ratediagram_simplified}, where $\tilde{I}$ includes both exposed and infectious states and $\tilde{R}$ includes the recovered, hospitalized, and dead states.  The model is further reduced to two dimensions by setting $\tilde{R} = N - \tilde{S} - \tilde{I}$.   When $N$ is large it is necessary to apply a course-grained discretization to the remaining state variables, $\tilde{S}$ and $\tilde{I}$.  We approximate the parameters for the simplified model so that the solution can be fed back to the full stochastic SEIHRD model.  We note that the suppression strategy could probably be approximated in one dimension with the assumption that $S=N$, and using a very high cost set at an arbitrary $\tilde{I}$ threshold.  Since, we aim to demonstrate that the strategy is the global optimum as compared also to mitigation strategies, we do not take this approach.

\begin{figure}
\centering
\includegraphics[width=0.5\textwidth]{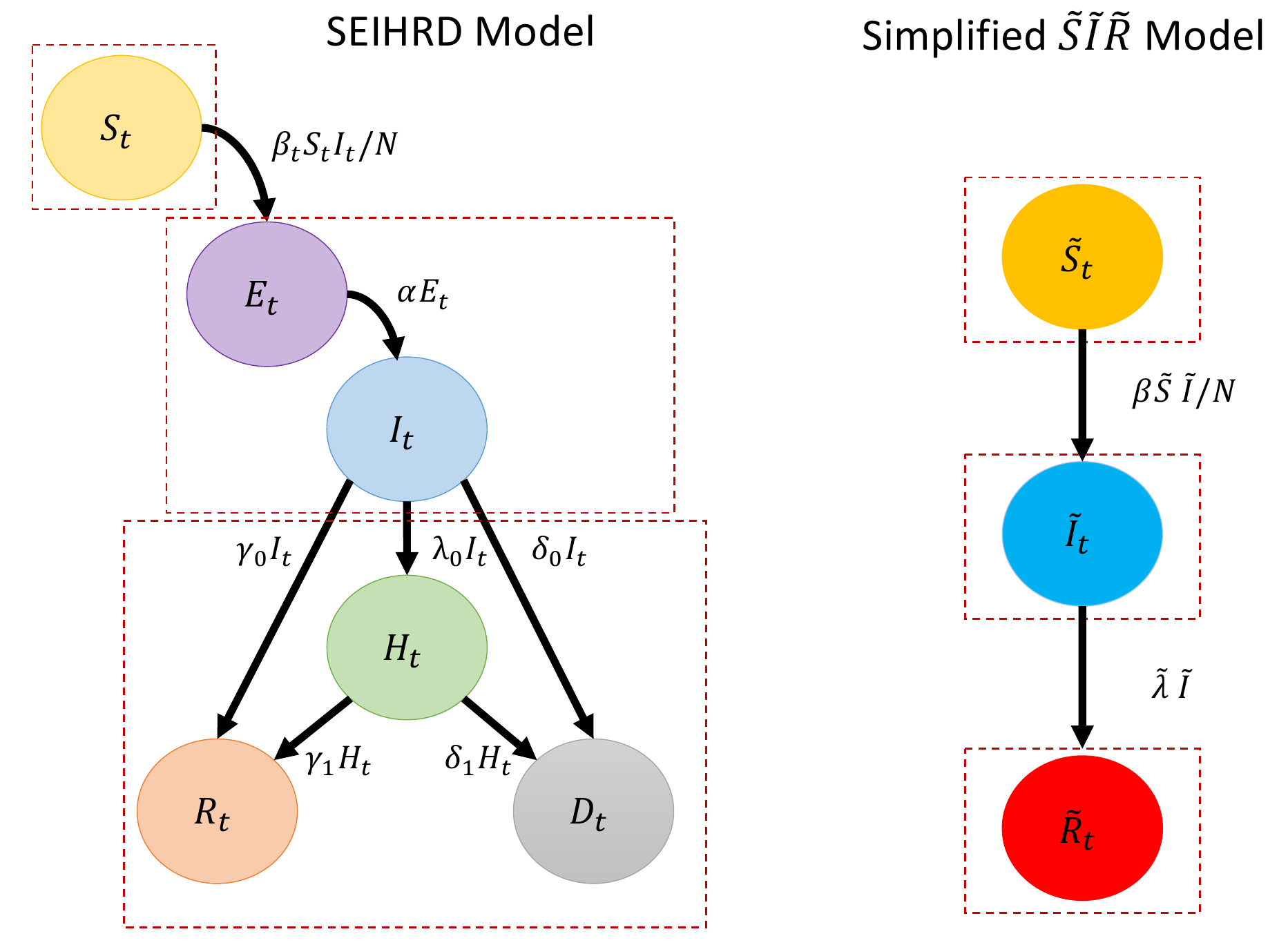}
\caption{Rate transition diagram for SEIHRD alongside the simplified $\tilde{S}\tilde{I}\tilde{R}$ model.}
\label{fig:ratediagram_simplified}
\end{figure}

For the simplified model there are now two states $\tilde{S}$ and $\tilde{I}$.  A transition occurs from $(\tilde{S},\tilde{I})$ to $(\tilde{S}-1, \tilde{I}+1)$ at rate $\beta\, \tilde{S}\, \frac{\tilde{I}}{N}$ (neglecting the effect of the incubation period), and a transition occurs from $(\tilde{S},\tilde{I})$ to $(\tilde{S}, \tilde{I}-1)$ at rate $\tilde{\lambda}\, \tilde{I}$, where $\tilde{\lambda} = \lambda_0+\gamma_0+\delta_0 \approx 0.217$.  We can then approximate $\tilde{H} = p_H\, \tilde{I}$ with $p_H = \frac{\lambda_0}{\gamma_1+\delta_1}\approx 0.071$ by considering the quasi-equilibrium when $\tilde{I}$ is near constant, and similarly $\tilde{D} = p_D\, \tilde{R}$ with $p_D=\frac{\delta_0 + \lambda_0 \left(\frac{\delta_1}{\gamma_1+\delta_1}\right)}{\tilde{\lambda}} 
%=\frac{\delta_0 +\frac{\delta_1}{\delta_1 + \gamma_1}\lambda_0}{\tilde{\lambda}}
\approx 0.0051$.

We run the model until the infection dies out, that is, until the first time $t$ is reached such that $\tilde{I_t}=0$.
%when $\tilde{I}=0$, $T = \inf\{t; \ I_t = 0\}$. 
We now assume $\beta$ has a feedback form, $$(\beta)=(\beta(\tilde{S},\tilde{I}))_{(\tilde{S},\tilde{I})\in \{0,\ldots, N\}^2}.$$   
  The cost is
$$
    J\big[(\beta)\big] = \mathbb{E}\Big[\int_0^TL\Big(\big(\beta(\tilde{S}_t,\tilde{I}_t)\big) + F(\tilde{H}_t)\Big)dt + G(\tilde{D}_T)\Big],
$$
where $\tilde{H}$ and $\tilde{D}$ are approximated as above, and $L$, $F$, and $G$ are the same as (\ref{eqn:costL}), (\ref{eqn:costF}), and (\ref{eqn:costG}).

We solve for the value function
$$
    V(\tilde{S}, \tilde{I}) = \sup_{(\beta)} \Big\{-\mathbb{E}\Big[\int_0^TL\Big(\big(\beta(\tilde{S}_t,\tilde{I}_t)\big) + F(\tilde{H}_t)\Big)dt + G(\tilde{D}_T)\Big];\ (\tilde{S}_0,\tilde{I}_0)=(\tilde{S},\tilde{I})\Big\},
$$
which satisfies $V(\tilde{S},0)=- G(\tilde{D})$ and solves the Bellman equation
$$
    \max_\beta\Big\{ \beta\, \tilde{S}\, \frac{\tilde{I}}{N}\big(V(\tilde{S}+1,\tilde{I}-1) - V(\tilde{S},\tilde{I})\big)+ \tilde{\gamma}\, \tilde{I}\, \big(V(\tilde{S}, \tilde{I}-1) - V(\tilde{S},\tilde{I})\big)-L(\beta)-F(\tilde{H})\Big\} = 0.
$$

It is an interesting challenge how to approximate the solution at a coarser discretization of the population variables in order to handle large $N$.  Special care must be taken near $\tilde{I}=0$ to `renormalize' the coefficients and accurately take into account the logarithmic behavior of the end time.  More details are given in Appendix \ref{apx:stochastic_numerics}.

\subsection{Stochastic Simplified Problem Results}

We plot five simulations of the SEIHRD model under `optimal control' from the simplified $\tilde{S}\tilde{I}\tilde{R}$ model (the correspondence between the simulations and the value function is only approximate), see Figure~\ref{fig:simopt} as a comparison to the deterministic suppression strategy in Figure~\ref{fig:optimal_mitigation}(a). As expected there is variability in the cost due to the fluctuations of the end time. For this approach we find only the globally optimal, `suppression', strategy and while the  exact values differ due to the coarse approximation, qualitatively the strategy is the same.
% To illustrate the `mitigation' strategy we plot the same simulation but with $d=\, $\$100,000 (maybe unreasonably low) in Figure \ref{fig:simopt_mitigation}.  The low death coefficient cost significantly reduces the cost of the `mitigation' strategy but has little effect on the cost of `suppression' strategy.  While there is still variability in the end-time with the `mitigation' solution, it no longer has much effect on the cost.  

\begin{figure}
$\includegraphics[width=\textwidth]{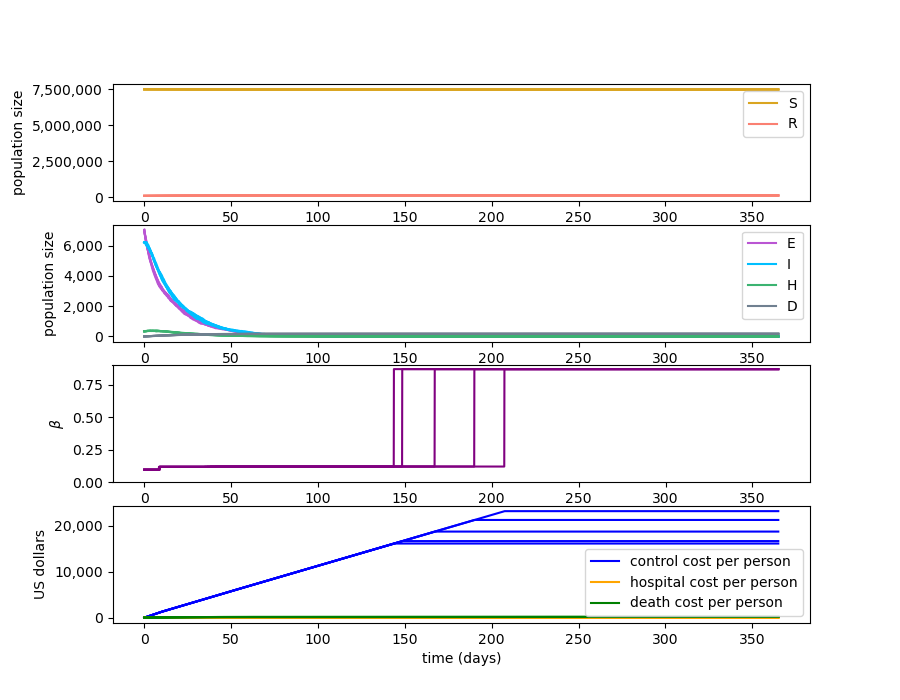}$\\

\caption{Five simulated optimal solutions, with the same parameters as in Figure~\ref{fig:optimal_mitigation}(a).  Note the fluctuations in the end time, while the other fluctuations are not noticeable. Individual simulations can be identified by matching the end time with the point the cost flattens.}
\label{fig:simopt}
\end{figure}

\begin{figure}
\begin{subfigure}[b]{0.48\textwidth}
$\includegraphics[width=\textwidth]{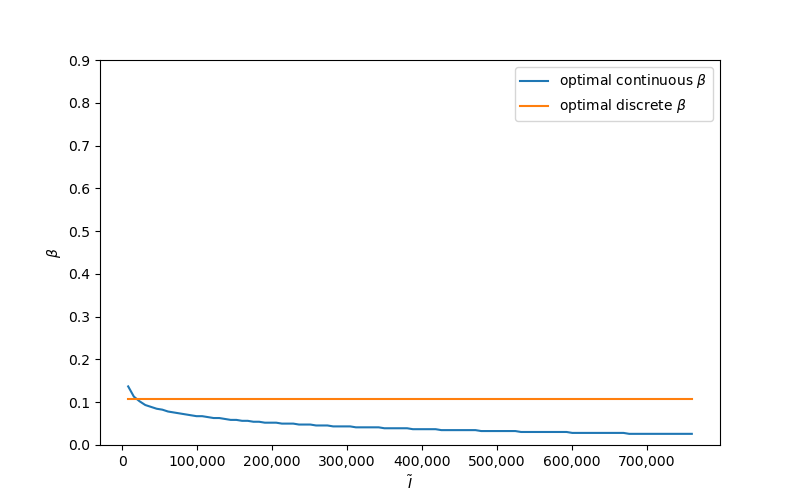}$\ 
\caption{Control $\beta$ versus infected population sizes with the susceptible population fixed.}
\end{subfigure}
\begin{subfigure}[b]{0.48\textwidth}
$\includegraphics[width=\textwidth]{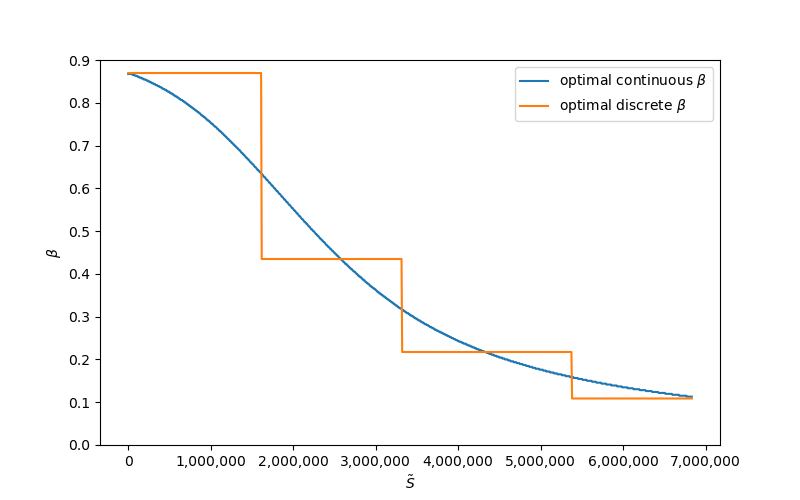}$
\caption{Control $\beta$ versus susceptible population sizes with the infected population fixed.}
\end{subfigure}
\caption{Continuous and discretized optimal $\beta$ values plotted as a function of  infected population sizes with the susceptible population fixed on the left (a), and as a function of susceptible population values with the infected population size fixed on the right (b). Same parameters as Figure \ref{fig:simopt}.}
\label{fig:beta_threshold}
\end{figure}

The main qualitative difference in the solutions from the finite time horizon, deterministic dynamics problem is that there is now a finite time where the infection dies out and the control returns to $\beta=b$.  

\subsection{Switching Times for Discrete $\beta$}

Instead of allowing continuous values for $\beta$ we now consider when $\beta$ is restricted to only four values calculated from ${\bf R_0} = 0.5, 1, 2$, and $4$. The motivation for discretizing $\beta$ is to reflect four policies that impact ${\bf R_0}$ through non-pharmaceutical interventions. 

The feedback control $\beta$, for both the continuous version and the discrete version, is plotted as a function of the infected population in Figure \ref{fig:beta_threshold}, where `suppression' and `mitigation' regimes can be identified.
 Indeed, the optimal discrete policy is very close to the policy with continuous $\beta$ rounded to the nearest admissible value. 
 We also plot the optimal $\beta$ for differing values of $\tilde{I}$ and $\tilde{S}$ in Figure \ref{fig:beta_threshold}(a) and (b), respectively.

The 
 ``switching time" between policies is shown in terms of infected and susceptible population size, providing feedback information on when to change policies.  For example, given a value for $\tilde{S}$, a decision-maker may identify the discretized level of optimal $\beta$.  As $\tilde{S}$ changes, Figure~\ref{fig:beta_threshold}(b) shows how this affects the feedback control policy, and when to switch to a different control level. 
 % In particular, the suppression strategy maintains $\beta$ at the lowest possible value until the disease has gone extinct, whereas the mitigation strategy with reduced death cost coefficient $d=\, $100,000 only reduces $R_0$ to $2$ at the peak of the epidemic.

\subsection{Vaccinations}

% { As a proof of concept, we introduce the possibility of vaccinations in the stochastic model. We include uncertainty in when vaccinations become available} by adding an additional discrete variable $u$ with binary values, $\{0,1\}$.  Here, $u=0$ represents that the vaccine is under development and not available, while 
% $u=1$ represents that the vaccine is being dispensed to the susceptible population at rate $o_1 \, N$.  We add this to the $\tilde{S}$ dynamics so that when $u = 1$, $(\tilde{S},\tilde{I}, 1)$ also transitions to $(\tilde{S} - 1, \tilde{I}, 1)$ at rate $o_1\, N$.  
% We assume that $(\tilde{S},\tilde{I},0)$ transitions to $(\tilde{S}, \tilde{I},1)$ at rate $w_0$. (Treating the time till vaccinations become available as an exponential random variable seems to be reasonable since there is so much uncertainty).  
We also demonstrate the stochastic model with vaccinations.  A new transition is added to reflect individuals leaving $\tilde{S}$ at a rate of $o\, N$.
The Bellman equation becomes
\begin{align}\label{eqn:bellman}
   0=&\ \sup_\beta\Big\{ \beta\, \tilde{S}\, \frac{\tilde{I}}{N}\big(V(\tilde{S}+1,\tilde{I}-1) - V(\tilde{S},\tilde{I})\big)\\
   &\ + \tilde{\gamma}\, \tilde{I}\, \big(V(\tilde{S}, \tilde{I}-1) - V(\tilde{S},\tilde{I})\big) \nonumber \\
    &\ + o\, N\, \big(V(\tilde{S} -1,\tilde{I})-V(\tilde{S},\tilde{I})\big)-L(\beta) - F(\tilde{H})\Big\}.\nonumber
\end{align}

%\subsection{Vaccination Results}

We solve the stochastic problem with vaccinations and run simulations using the same parameters from the US COVID-19 data as of January 1, 2021, (\ref{eqn:initial_values_US}).
As shown in the top graph in Figure \ref{fig:vacc_sim}, the susceptible population $S$ drops off linearly to near zero, as the vaccination takes affect. Comparing Figure~\ref{fig:simopt} with Figure~\ref{fig:vacc_sim}, it is clear that the solutions are qualitatively the same, as the control $\beta$ starts off low and increases later, while  cost is largely due to the control cost in both situations.

With the same parameters we view how the vaccination affects the $\beta$ thresholds in Figure \ref{fig:vacc_d9}.  The overall solution remains close to the solution without vaccinations of Figure \ref{fig:simopt}.
Figure \ref{fig:vacc_d9} plots the effect on the threshold of $\beta$ values, with and without a vaccine. A surprising feature is that the optimal $\beta$ is almost always lower with a vaccine compared to without one.  While vaccination allows for a faster relaxation of the control, it is only due to the decrease of the susceptible population and at the same susceptible population a stricter control should be used alongside the vaccinations.

\begin{figure}
$\includegraphics[scale=0.57]{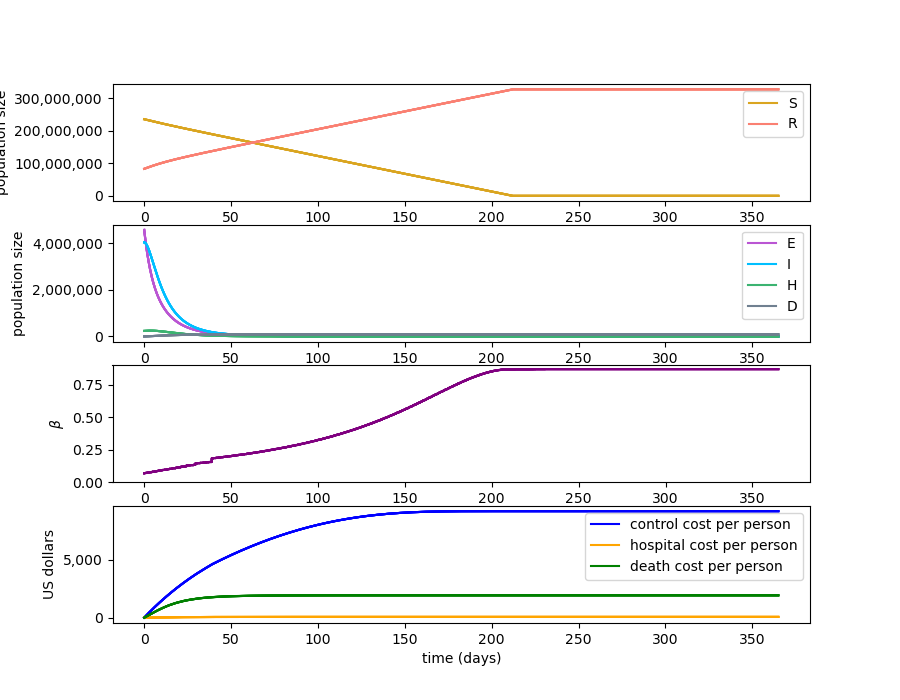}$
\caption{Simulations for the stochastic model with vaccinations.  Parameters are the same as Table~\ref{table:parameters} and Table~\ref{table:control_parameters} with initial conditions from (\ref{eqn:initial_values_US}) reflecting the US COVID-19 data as of January 1, 2021, and $o=1/300$. 
\label{fig:vacc_sim}}
\end{figure}

\begin{figure}
\begin{subfigure}[b]{0.48\textwidth}
$\includegraphics[width=\textwidth]{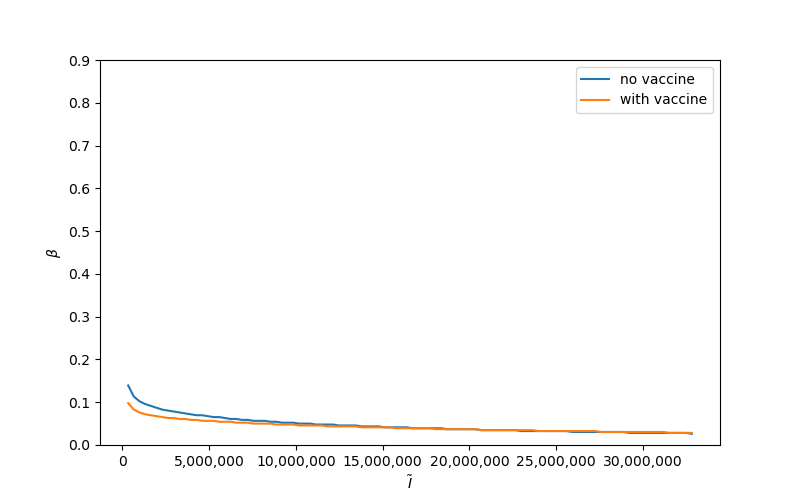}$\ 
\caption{Control $\beta$, with and without vaccine, versus infected population sizes with the susceptible population fixed.}
\end{subfigure}
\begin{subfigure}[b]{0.48\textwidth}
$\includegraphics[width=\textwidth]{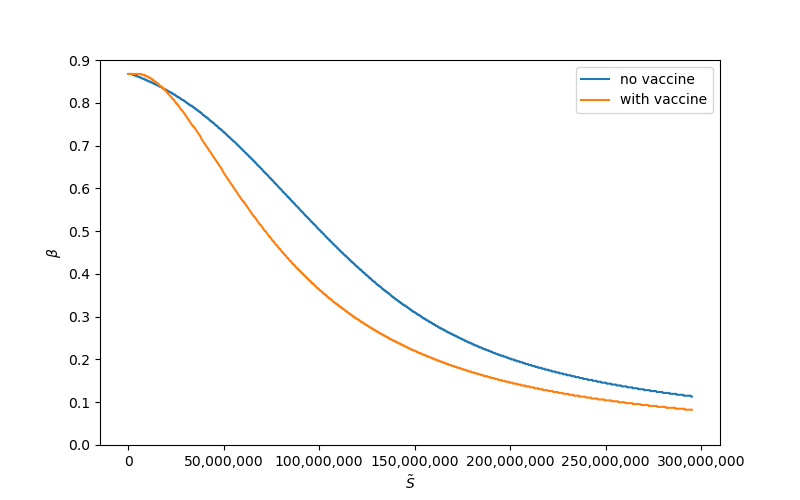}$
\caption{Control $\beta$, with and without vaccine, versus susceptible population sizes with the infected population fixed.}
\end{subfigure}
\caption{Optimal $\beta$ thresholds, with and without a vaccine, with the susceptible population fixed at $\tilde{S} =\, $294,000,000 in (a), and with the infected population fixed at $\tilde{I}=\, $6,500,000 in (b).  Parameters are the same as Table~\ref{table:parameters} and Table~\ref{table:control_parameters} with initial conditions from (\ref{eqn:initial_values_US}) reflecting the US COVID-19 data as of January 1, 2021, and $o=1/300$. 
\label{fig:vacc_d9}}
\end{figure}

   \subsection{Stochastic Control Numerics}\label{apx:stochastic_numerics}
  
        The Bellman equations (\ref{eqn:bellman}) can be solved in a single sweep of value iterations making sure that we first increase $\tilde{S}$ then increase $\tilde{I}$, and with vaccinations we start with $u=1$ and then do $u=0$.  We discretize $\beta$ in increments of $0.01$.  

        When $N$ is much larger than 1000 we do not solve the equations directly as the computational cost is of order $N^2$.  Instead we discretize the population variable in 1000 increments.  We let $\Delta k = 0.001 \, N$ denote the discretization increment (which we assume is greater than 1).  It is then possible to solve approximate Bellman equations, where we replace 
        $$
            V(\tilde{S},\tilde{I}-1,u) - V(\tilde{S},\tilde{I},u) \approx \frac{V(\tilde{S},\tilde{I}-\Delta k,u) - V(\tilde{S},\tilde{I},u)}{\Delta k}.
        $$

        Unfortunately, this leads to a bad approximation when $\tilde{I}$ is small.  For example the expected time to transition from $\tilde{I} = \Delta k$ to $\tilde{I} = 0$ with $\tilde{S}=0$ is approximated by $\frac{1}{\alpha}$, whereas the correct expected time can be computed as
        $$
            \frac{1}{\alpha}\sum_{j=1}^{\Delta k} \frac{1}{j},
        $$ 
        which is about 7.5 times larger when $\Delta k=1000$. Since $\tilde{I}=\Delta k$ must be visited by any solution before the end time, this error would propagate through the whole problem.
         Because of this we renormalize the coefficient $\alpha$ at state $\tilde{I}$ using the formula 
        $$
            \frac{1}{\alpha} \rightarrow \frac{1}{\alpha} \Big(\sum_{j=\tilde{I} - \Delta k +1}^{\tilde{I}} \frac{1}{j}\Big).
        $$
        The sum is approximated using the standard formula
        $$
            \sum_{j=1}^{k} \frac{1}{j} \approx \log(k)+\gamma_e+\frac{1}{2k}-\frac{1}{12k^2}
        $$
where $\gamma_e \approx 0.577$ is the Euler-Mascheroni constant.

\bibliographystyle{siamplain}
\bibliography{SEIRBib}

\end{document}